\documentclass[notitlepage,leqno,11pt]{article}
\usepackage{latexsym}
\usepackage{amsmath}
\usepackage{amsfonts}
\usepackage{amssymb}
\usepackage{amscd}
\renewcommand{\a }{\alpha }
\renewcommand{\b }{\beta }
\renewcommand{\d}{\delta }
\newcommand{\D }{\Delta }
\newcommand{\e }{\varepsilon }
\newcommand{\g }{\gamma}
\renewcommand{\i }{\iota}

\renewcommand{\l }{\lambda }
\renewcommand{\L }{\Lambda }

\newcommand{\n }{\nabla }

\newcommand{\s }{\sigma }
\newcommand{\Sig }{\Sigma}

\renewcommand{\o }{\omega }

\newcommand{\ov}{\overline}
\newcommand{\intbar}{\mathop{\int\makebox(-13.5,0){\rule[4pt]{.7em}{0.3pt}}%
\kern-6pt}\nolimits}
\newcommand{\wtilde }{\widetilde}

\newenvironment{pf}{\noindent{\sc Proof}.\enspace}{\rule{2mm}{2mm}\medskip}
\newenvironment{pfn}{\noindent{\sc \bf Proof }}{\rule{2mm}{2mm}\medskip}

\newtheorem{remark}{Remark}[section]
\newcommand{\R}{\mathbb{R}}

\newcommand{\N}{\mathbb{N}}
\renewcommand{\o }{\omega }

\author{ Mohameden  Ahmedou \; $\&$\;\; Mohamed Ben Ayed }

\date{}

\title{\bf  Morse inequalities at infinity for a resonant mean field  equation}
\begin{document}

\newtheorem{lem}{Lemma}[section]
\newtheorem{pro}[lem]{Proposition}
\newtheorem{thm}[lem]{Theorem}
\newtheorem{rem}[lem]{Remark}
\newtheorem{cor}[lem]{Corollary}
\newtheorem{df}[lem]{Definition}

\maketitle

\centerline{   \emph{In fond memory of Abbas Bahri }}

\begin{center}
{\bf Abstract}
\end{center}
In this paper we study the following  mean field type  equation
\begin{equation*}
(MF) \qquad  -\D_g u \, =  \varrho ( \frac{K e^{u}}{\int_{\Sig} K e^{u} dV_g}  \, - \, 1) \, \mbox{ in } \Sigma,
\end{equation*}
 where $(\Sigma, g)$ is a closed oriented surface of unit volume $Vol_g(\Sigma)$ = 1,  $K$ positive smooth function and $\varrho= 8 \pi m$,  $ m \in \N$.
 Building on the critical points at infinity approach initiated  in  \cite{ABL17} we   develop, under generic condition on the function $K$ and the metric $g$,  a full Morse theory by proving Morse inequalities relating the Morse indices of the critical points, the indices of the critical points at infinity, and the Betti numbers of the space of formal barycenters $B_m(\Sigma)$.\\
We derive from these  \emph{Morse inequalities at infinity} various new  existence  as well as multiplicity results  of the mean field equation in the resonant case, i.e. $\varrho \in 8 \pi \N$.
\begin{center}

\bigskip
\noindent{\bf Key Words:}  Critical points at infinity, Morse theory, Space of formal  barycenters.

\centerline{\bf AMS subject classification:  35C60, 58J60, 35J91.}

\end{center}

%%%%%%%%%%%%%%%%%%%%%%%%%%%%%%%%
\section{Introduction and statement of the results}
%%%%%%%%%%%%%%%%%%%%%%%%%%%%%%%%

Let $(\Sigma, g)$ be a closed and  oriented surface of unit volume $Vol_g(\Sigma)$ = 1. For  $\varrho\in \R^+$ real number and $K$ positive smooth function, we consider the following Mean Field equation
\begin{equation}\label{eq:TS}
(MF) \qquad -\D_g u \, =  \varrho ( \frac{K e^{u}}{\int_{\Sig} K e^{u} dV_g}  \, - \, 1) \, \mbox{ in } \Sigma .  \end{equation}

Problem $(MF)$ arises as limiting equation in some mean field approximation in the study of limit point vortices of Euler flows, in Onsager vortex theory, and also in the description of self dual condensates in some Chern-Simons-Higgs models. See for example the papers  \cite{AP15, clmp1, clmp2, cl1, cl2, djlw, k, YYan, Lei1, st}, and the monographs of Tarantello \cite{Tbook} and of Yang \cite{Ybook} as well as  the references therein. Furthermore  its study  is also  motivated by the prescribed Gauss curvature problem  in differential geometry. Indeed, if  $\Sig$ is the standard $2$-sphere  $\mathbb{S}^2$ endowed with its standard metric $g_{stand}$ and $\varrho = 8 \pi$,   then Problem $(MF)$ is the so called \emph{ Nirenberg's  problem } in conformal geometry, which consists of finding a metric $g$ conformally equivalent to $g_{stand}$ and  whose Gauss curvature is given by the function $K$. See  \cite{CY, Han, CGY} and the references therein.\\
Equation $(MF)$ has a variational structure and the associated  Euler Lagrange functional   $J_{\varrho}$ is defined by
$$ J_{\varrho}(u) \, := \, \frac{1}{2} || u ||^2 \, - \varrho \ln( \int_{\Sigma} K e^{u} dV_g) \quad \mbox{ for }  u \in \mathring{H ^1 }(\Sigma)$$
where
$$ ||u|| :=    ||\n u||_{L^2}  $$
is the norm used to equip
\begin{equation}\label{eq:H10}
\mathring{H}^1(\Sigma) := \Big \{   u \in H^1(\Sigma): \int_{\Sigma} u dV_g \, = 0\Big\}.
\end{equation}
It turns out that the analytical aspect of this variational problem depends on the range of values taken by the parameter $\varrho$. Indeed
\begin{itemize}
 \item If $\varrho \notin 8 \pi \mathbb{N}$ then the associated variational problem is compact in the sense that changes of the topology of the  sublevel sets of $J_{\varrho}$ are induced only by critical points. We call it the \emph{non-resonant case.}
  \item If $\varrho \in  8 \pi \mathbb{N}$ then the associated variational problem is not compact, i.e.\ changes in the difference of topology between its sublevel sets might be induced by critical points at infinity, these are non compact orbits of the gradient flow whose level sets are bounded (see \cite{bah} for a related notion for Yamabe type problem).  We call it the \emph{resonant case.}
\end{itemize}
In the non-resonant case, Yanyan Li \cite{YYan} proved that the Leray-Schauder degree of the solutions of $(MF)$ is well defined and is a topological invariant. He also advised a method to compute it by analyzing the \emph{ degree  jump } across the critical values $8 m \pi$. C.C. Chen and C.S. Lin \cite{cl1, cl2}, proved a priori estimates for the solutions,  and used the strategy advised by Yanyan Li to compute the Leray-Schauder degree. Then they derived that Problem $(MF)$ is solvable provided that the surface has a positive genus. Later Z. Djadli  \cite{dj}, using the topological argument of \cite{DM-annals}, proved the existence of a solution of $(MF)$ in this case for surfaces of  every genus. Furthermore  A. Malchiodi \cite{andrea} developed  a Morse theory in this case  and F. De Marchis \cite{fdM} proved full Morse inequalities as well as multiplicity results for generic data $(K,g)$.\\
For the \emph{resonant case} the main difficulty lies in the fact that to find critical points of $J_{\varrho}$  one has to understand the contributions of the critical points at infinity to the topology of its sublevel sets. This amounts to proving a Morse lemma at infinity in a highly degenerate setting where gradient flow lines may converge only in the sense of measures to a sum  of Dirac masses. Based on ideas going back to Lions \cite{Lions:1985a} and Struwe \cite{Struwe:1984} one can often perform a "Lyapunov-Schmidt reduction at infinity", and finds that the Dirac deltas are located at critical points of a finite-dimensional reduced function. Such an approach has been developed by the authors in collaboration with M. Lucia in \cite{ABL17}. There,  a Morse reduction in the neighborhood of the critical points has been performed and the induced difference of topology between the levels of $J_{\varrho}$ has been computed. In this paper the authors want to develop further  this approach by proving \emph{Morse inequalities} between \emph{critical points at Infinity} and use them to deduce new existence and multiplicity results.\\
To state our main results  we introduce the following notation and assumptions:\\
For $m \in \N  $,  we let
$$\mathbb{F}_m(\Sigma):= \left  \{  (a_1, \cdots,a_m) \in \Sigma^m; a_i \neq a_j \mbox{ if } i \neq j \right  \} $$
denote the space of configurations  of $\Sigma^m$  and we define the following reduced energy functional
$$ \mathcal{F}^{K}_m: \, \mathbb{F}_m(\Sigma) \, \,  \longrightarrow \,    \R, \quad \mbox{ by }$$
\begin{align}  & \mathcal{F}^{K}_1(a) \, := \,  \ln K(a) \, + 4 \pi  H(a,a)  \quad (\mbox{for }m=1; \,  a \in \Sigma),  \label{eq:F1}\\
  & \mathcal{F}^{K}_m(A) \, := \,  \sum_{i=1}^{m}  \bigg( \ln K (a_i) \, + \, 4 \pi  H(a_i,a_i) \, + 4 \pi \sum_{j \neq i} G(a_i,a_j)  \bigg)\,  (\mbox{for }  m \geq 2 ;\,  A = (a_1,\cdot,a_m)), \label{eq:Fm}
\end{align}
where $G$ is the Green's function of $ -\D_g$ and  $H$ its regular part. See \eqref{eq:Green} and \eqref{eq:GH} for the precise definition.
We notice  that the function $  \mathcal{F}^{K}_m$ achieves its infimum  for each $m\geq 1$ and for $m=1$ it also  achieves its maximum.

Next for $A =(a_1,\cdots,a_m) \in \mathbb{F}_m$ a critical point of $ \mathcal{F}^{K}_m$ we define  a vector $(\mathcal{F}^A_1, \cdots, \mathcal{F}^A_m)$, where
$$ \mathcal{F}^A_i: \, \Sigma  \longrightarrow  \R, $$
are real valued functions defined as follows:
\begin{align}
&  \mathcal{F}_1^a(x) \, := \, K(x) \exp(  \, 8 \pi H(a,x) ) , \quad \mbox { if } m = 1,  \label{f1} \\
& \mathcal{F}^A_i(x) \, := \, K(x) \exp\Big( \, 8 \pi H(a_i,x) \, + \, 8 \pi \sum_{j \neq i} G(x,a_j)\Big) \quad \mbox{ if } m\geq 2 \, \mbox{ for } \, i=1, \cdots,m .\label{eq:fm}
\end{align}
Moreover we associate to every critical point $A =(a_1,\cdots,a_m)$ of $ \mathcal{F}^{K}_m$ the following quantity
\begin{equation}\label{eq:LA}
  \mathcal{L}(A) \, := \, \sum_{i=1}^{m} \left ( \D_g \mathcal{F}^A_i(a_i) - 2 K_g(a_i) \, \mathcal{F}^A_i(a_i) \right ),
\end{equation}
where $K_g$ denotes the Gauss curvature of $(\Sigma, g)$.
Next we define the following set
\begin{equation}\label{eq:K-}
\mathcal{K}^{-}_m \, := \, \left \{   A \in \mathbb{F}_m(\Sigma): \, \n \mathcal{F}^{K}_m(A) \, = \, 0, \, \mathcal{L}(A) < 0 \right \}.
\end{equation}
%%%%%%%%%%%%%%%%%%%%%%%%%%%%%%%%%%%%%%%%%%%%%%%%%%%%%%%%%
 Our first type of results are based on the construction of solutions of sub- and sup-approximation with fixed morse indices, which we  show that, under appropriate assumptions of critical points of the related reduced energy, they  give rise to solutions of Equation $(MF)$.

\begin{thm}\label{th:crit1}  Let $m\geq 1$  and $\varrho = 8 \pi m$ and assume that the reduced energy  $\mathcal{F}^K_m$ has only non degenerate critical points. Then
\begin{enumerate}
  \item[1)]
  if at every local minimum $A \in \mathbb{F}_m(\Sigma)$  of $\mathcal{F}^K_m$ we have that $\mathcal{L}(A) < 0.$ Then Equation $(MF)$ has at least one solution whose generalized Morse index is $ 3m$,
  \item[2)] for $m \geq 2$,
  if at every critical point  $A \in \mathbb{F}_m(\Sigma)$  of $\mathcal{F}^K_m$ with  morse index $2$ we have that $\mathcal{L}(A) > 0$. Then Equation $(MF)$ has at least one solution whose generalized Morse index is $ 3m-3$,
\end{enumerate}
where  the \emph{generalized Morse index}  of a solution $\o$ of $(MF)$ is the  sum of its  Morse index and the dimension of the Kernel of the linearized operator
$$
\mathcal{T}_{\o}(\varphi):=  - \D_g \varphi \, - \, \varrho \frac{K e^{\o} \varphi}{\int_{\Sigma} K e^{\o}}.
$$
\end{thm}

We observe  that, for $m=1$,  the functional $J_{8\pi}$ is bounded below and therefore if there exists a local maximum point $A$ of $\mathcal{F}^K_1$  such that $\mathcal{L}(A) > 0$ then $J_{8\pi}$ achieves its minimum (see \cite{ABL17}, Theorem 1.1).  We point out that there is no such points for the Nirenberg's problem on $\mathbb{S}^2$ and more generally there are no such points if the surface has a positive Gauss curvature.

\noindent
Our next and main result of this paper is to prove \emph{Morse inequalities}. To that aim we  introduce the following non degeneracy assumption:

\noindent \emph{  We say that the pair  $(g, K)$ satisfies the condition $(\mathcal{N}_m)$ if the following conditions are satisfied:
\begin{enumerate}
\item[(i)]  The critical points of $\mathcal{F}^{K}_m$ are non degenerate and for  every critical point $A$, we have $\mathcal{L}(A) \neq 0,$
  \item[(ii)] All the critical points of $J_{8 \pi m }$ are non degenerate.
\end{enumerate}  }

\begin{remark}
The non degeneracy condition $(\mathcal{N}_m)$ is generic. Indeed, denoting by $\mathcal{M}$ the set of riemannian metrics on $\Sigma$ and $ C^{0,1}_+(\Sigma)$ the space of positive Lipschitz functions on $\Sigma$. It follows  from   transversality arguments, there is an open and dense subset $\mathcal{D} \in \mathcal{M} \times C^{0,1}_+$ such that for $(g,h) \in \mathcal{D}$, the functional $J_{\varrho}$ is a Morse function. See for example \cite{fdM}, pp. 2179-80, for an argument based on some generic properties for nonlinear boundary value problems proved by Saut and Temam \cite{STEM}. Moreover the condition $(i)$ is also a generic one.
\end{remark}
Furthermore, under the non degeneracy condition $(i)$ of $(\mathcal{N}_m)$,  we associate  to every $A \in \mathcal{K}^-_m$  an index
\begin{equation}\label{eq_index} \iota_{\infty}:\mathcal{K}^- _m\to \N \quad \mbox{ defined  by } \quad
  \iota_{\infty}(A) \, :=\, 3 m -1 - Morse(\mathcal{F}^{K}_m, A),
\end{equation}
where $ Morse(\mathcal{F}^{K}_m, A)$ stands for the Morse index of $\mathcal{F}^{K}_m $ at its critical point $A$.

\noindent
Next for $m \in \N$, let $$B_m(\Sig):= \left \{ \sum_{i=1}^{m} \gamma_i \d_{a_i}, \, a_i \in \Sigma, \, \gamma_i \geq 0; \, \sum_{i=1}^{m} \g_i\,  = \, m \right \}  $$
denote the set of formal barycenters of order $m$. For $i\geq 0$, we set
$$ \beta^m_i:= \, rank \, \, H_i( B_m(\Sigma), \mathbb{Z}_2). $$
Then  $B_m(\Sig)$ is a stratified set  of top dimension $3m-1$ and therefore $\beta_i^m = 0$ for $i \geq 3m$.

\noindent
In the next result we provide Morse inequalities relating the Morse indices of the solutions of $(MF)$, the $\iota_{\infty}$-indices of the critical points at infinity and the Betti numbers  $\beta_i^m$ of the set of formal barycenters $B_m(\Sig )$. Before stating these inequalities, we recall that it follows from the blow up analysis of solutions of $(MF)$, see \cite{LS,cl1, cl2} that their  Dirichlet energy is uniformly bounded. See  please  statement $(ii)$ in  Theorem 3.1, p. 1686 in \cite{cl2}. Hence it follows from Moser-Trudinger inequality and elliptic theory that solutions of $(MF)$ with zero mean value integral are uniformly bounded and hence their Morse indices are uniformly bounded.\\
We are now ready to state our Morse inequalities in the case $ \varrho= 8 \pi m; m \geq 2$:
\begin{thm}\label{th:supc3}
Let $\varrho= 8 \pi m$, $ m \geq 2$ and assume that the function $K$ satisfies the non degeneracy condition $(\mathcal{N} _m)$ and let $\ov{m}$ be the highest Morse index of the critical points of  $J_{\varrho}$. Then the following Morse inequalities hold
\begin{enumerate}
  \item[(a)]
  For every $ 2 \leq k \leq \ov{N}:=\max (3m -1  , \ov{m})$ there holds
  \begin{equation}\label{azertnew}   \nu_k  \, + \nu^{\infty}_k \, \geq  \,   \beta_{k-1}^{m-1}, \end{equation}
  where $\nu_i$ denotes the number of critical  points of $J_{\varrho}$ with Morse index $i$ and
   \begin{equation} \label{mq} \nu^{\infty}_q:= \#\{ A \in \mathcal{K}^-_m; \iota_{\infty}(A) = q\}. \end{equation}
  \item[(b)] Let $ \chi(\Sigma)$ be the Euler Characteristic of $\Sigma$, it holds
  $$
  \sum_{i=0}^{\ov{m}} (-1)^{i  } \nu_i  \, + \sum_{A \in \mathcal{K}^{-}_m} (-1)^{\iota_{\infty}(A) } \, = \, \binom{m -1- \chi(\Sigma)}{m-1}.
  $$
\end{enumerate}
In particular  the number of critical points of $J_{\varrho}$ is lower bounded by
$$
\left   |   \binom{m - 1 - \chi(\Sigma)}{m-1}   -  \sum_{A \in \mathcal{K}^{-}_m} (-1)^{\iota_{\infty}(A) }    \right |.
$$
\end{thm}

\noindent
As corollary of the  above \emph{ Morse inequalities at infinity } we derive the following existence result:

\begin{thm}\label{th:supc4}
Let $m\geq 2$ and $\varrho= 8 \pi m$   and assume that the function $K$ satisfies the non degeneracy condition $(\mathcal{N} _m)$. Suppose  there exists $2 \leq q_0 \leq 3m -3$ such that \begin{itemize}
  \item there is no element in $\mathcal{K}^{-}_{m}$ of index $q_0$
  \item $\b_{q_0-1}^{m-1} \neq 0$.
\end{itemize}
Then Equation $(MF)$ has at least one solution whose Morse index is $q_0$.\\
In particular, since $\beta_{3m-4}^{m-1}= H_{3m-4}(B_{m-1}(\Sigma)) \neq 0$, if there is no   element in $\mathcal{K}^{-}_{m}$ of index $3m-3$, then Equation $(MF)$ has at least one solution whose Morse index is $3m-3$. Furthermore, the number of solutions is lower bounded by $\beta_{3m-4}^{m-1}$.
\end{thm}
As a complement of the statement made  in the above theorem, we prove  the following result:
\begin{thm}\label{th:new}
Let $m\geq 2$ and $\varrho= 8 \pi m$ and  assume that the function $K$ satisfies the non degeneracy condition $(\mathcal{N} _m)$.
Suppose there exists $q_0$ such that \begin{enumerate}
\item there exists $Q_0 \in \mathcal{K}^-_m$ with $\iota_\infty(Q_0) = q_0$,
\item for each $Q \in  \mathcal{K}^-_m$, we have $\iota_\infty(Q) \notin \{ q_0 -1, q_0 +1\}$,
\item $\displaystyle{ \beta_{q_0-1}^{m-1} \, = \, 0 } $.
\end{enumerate}
Then Equation $(MF)$ has at least one solution.
\end{thm}

\begin{remark}
We notice that  the indices of the  critical points at infinity are lower bounded by $m-1$.  Hence  for $m=4$ suppose
\begin{enumerate}
\item there exists $ Q_0 \in \mathcal{K}^-_m$ such that $ \iota_{\infty}(Q_0) = 3$ (i.e. $\mathcal{F}_m^K(Q_0)$ is a local  maximum of $ \mathcal{F}_m^K$),
\item for all $Q \in \mathcal{K}^-_m$, we have $ \iota_{\infty}(Q) \neq 4$,
\end{enumerate}
then it  follows from Theorem \ref{th:new} that Equation $(MF)$ has at least one solution.
\end{remark}
More generally any violation of the  above \emph{Morse inequalities at infinity} of Theorem \ref{th:supc3} induces the existence of a solution to Equation $(MF)$. In particular we have the following existence result:
\begin{thm}\label{th:supc5}
Let $m\geq 2$ and $\varrho= 8 \pi m$   and assume that the function $K$ satisfies the non degeneracy condition $(\mathcal{N} _m)$ . If  one of the following  inequalities is not satisfied
\begin{itemize}
  \item $\nu ^{\infty}_{3m-2} \geq \nu ^{\infty}_{3m-1}$.
  \item $\nu ^{\infty}_{3m-3} - \beta_{3m-4}^{m-1}\geq \nu ^{\infty}_{3m-2} - \nu ^{\infty}_{3m-1} $.
\end{itemize} Then Equation $(MF)$ has at least one solution (where $\nu_q^\infty$ is defined in \eqref{mq}).
\end{thm}
\noindent
The remainder of this paper is organized as follows: we set up some notation and definitions in Section 2, introduce an $\e-$neighborhood of potential  critical points at infinity in Section 3 and   provide in Section 4 useful expansions of the gradient of the Euler-Lagrange functional in \emph{ the neighborhood at infinity}. Section 5 is devoted to the characterization of such critical points at infinity and we provide the proof of our main results in Section 6. Finally we provide useful estimates of the projected bubble (see \eqref{eq:phial} for the definition) in the appendix.

%%%%%%%%%%%%%%%%%%%%%%%%%%
\section{Notation and definitions}
%%%%%%%%%%%%%%%%%%%%%%%%%%%

To state our results we need to fix some notation.\\
For  $\xi \in \Sigma$, let  $y_\xi$ be  a local chart defined in  a neighborhood of $\xi$ onto $B(0, 2\eta_0)$ such that  $y_\xi(\xi)=0$. In the sequel, we will denote $B_\xi (\eta):= y_\xi ^{-1} (B(0,\eta))$. Furthermore,
 thanks to  the isothermal coordinates, we have that for every $a \in \Sigma$, there exists a function $u_a \in C^{\infty}(\Sigma)$, satisfying $u_a(a) = 0$ and $\n u_a(a) = 0$, such that for the conformal metric $g_a := e^{u_a} g$, there hold
$$ g_a \, = \, \d_{ij} \mbox{ in } B_a(2\eta_0),\, \,  \D_g = e^{u_a} \D_{g_a}, \, \,  dV_{g} = e^{-u_a} dV_{g_a}, \, \mbox{ and } \,  \D_{g_a} u_a = 2K_g(.)e^{-u_a} \mbox{ in } B_a(2\eta_0), $$
where $0 < \eta_0 < \eta_a < \frac{1}{4} \, inj_{g_a}(\Sigma)$ where $inj$ stands for the injectivity radius.\\
\noindent
Next for  $a \in \Sigma$, let $G(a,.)$ be the Green's function of $\D_g$ defined by
\begin{equation}\label{eq:Green}
\begin{cases}
  - \D_g \, G(a,x) \, + \, 1    & = \, \d_a \,  \mbox{in }  \Sigma , \\
  \int_{\Sigma} G(a,x) dVg(x) & = \, 0.
\end{cases}
\end{equation}
For $0 < \eta < \frac{1}{4}  \eta_0$ and $a\in \Sigma$  we define the smooth  cut-off function
\begin{equation}\label{eq:cutof}
  \psi_a(x):= \psi(| y_a(x) | ) \quad \mbox { where } \quad \psi(t) :=
\begin{cases}
  t, & \mbox{if } t \in [0,{\eta}], \\
  2\eta, & \mbox{if } t \geq {2 \eta}, \\
  \psi(t) \in [\eta, 2\eta], \psi \in C^{\infty} & \mbox{otherwise}.
\end{cases}
\end{equation}
Next  for $\xi\in \Sigma$, as usual we decompose  $G(\xi,.)$ as follows
\begin{equation}\label{eq:GH} G(\xi,x)= -\frac{1}{2\pi} \ln ( \psi_\xi(x) ) + H(\xi,x).\end{equation} From the definition of $G(\xi,.)$ we derive that $H(\xi,.)$ has to satisfy
\begin{equation}\label{eq:H}
\D_g H(\xi,.) \, - 1 \,  =  \D_g G(\xi,.) - 1 + \frac{1}{2\pi} \D_g \ln(\psi_\xi(.))
= \begin{cases}  0 \mbox{ in } (\Sig\setminus B_\xi(2\eta))\cup B_\xi(\eta), \\
                            \frac{e^{u_\xi}}{2\pi} \D_{g_\xi} \ln(\psi_\xi(.)) \mbox{ in }  B_\xi(2\eta) \setminus B_\xi(\eta).\end{cases}
\end{equation}
%%%%%%%%%%%%%%%%%%%%%%%%%%%%%%%%%%%%%
In the next  section, we will describe the lack of compactness occurring in the variational problem associated to Equation $(MF)$.  To do so we need to introduce some highly concentrated functions, the so called \emph{bubbles} and a neighborhood of such bubbles.\\
First, we recall  that the following functions
$$ \wtilde{\d}_{a,\l}(x) := \ln( \frac{8 \l^2}{(1 + \l^2  |a - x |^{2})^2}),  \mbox{ for } a\in \R^2, \mbox{ and } \l > 0,$$
are the solutions of $$ -\D u = e^u \quad \mbox{in } \R^2 \quad \mbox{ with } \quad \int_{\R^2} e^u < \infty.$$

Next for $a \in \Sigma$ and $\l > 0$ we define  the standard bubble   as
\begin{equation}\label{eq:delta}
\d_{a,\l} (x):= \ln( \frac{8 \l^2}{(1 + \l^2   \psi_a(x) ^{2})^2}),
\end{equation}
where $\psi_a$ is defined in \eqref{eq:cutof}. In order to construct a suitable \emph{approximated solution} of $(MF)$ we follow the strategy introduced by A. Bahri and J.M. Coron in their study of the Yamabe equation \cite{BCd} and we  introduce the \emph{ projected bubble } $\varphi_{a,\l}$ to be the unique solution of
\begin{equation}\label{eq:phial}
\begin{cases}
-\D_g \varphi_{a,\l} \, = \,  e^{\d_{a,\l} + u_a}\, - \, \int_{\Sigma}  e^{\d_{a,\l} + u_a} dV_g  \,  \mbox{ in } \Sig ,\\
\int_{\Sig} \varphi_{a ,\l}\,  dV_g \, = \, 0.
\end{cases}
\end{equation}
We observe that for $A:=(a_1, \cdots,a_m) \in \Sig ^m$ and $\L:= (\l_1,\cdots,\l_m) \in (\mathbb{R}^+)^m$ the  sequence of functions $(U_{A,\L})_{\L}$ defined by
\begin{equation}\label{eq:Ual}
  U_{A,\Lambda}:= \sum_{i=1}^{m} \varphi_{a_i,\l_i}
\end{equation}
is a Palais-Smale sequence for Equation $(MF)$ for $\varrho= 8 \pi m$  if $A$ and $\L$ satisfy the following balancing condition
\begin{equation}\label{eq:balanc-cond}
\forall i \neq j \quad  \l_i^2  \mathcal{F}^A_i(a_i) \, = \, \l_j^2  \mathcal{F}^A_j(a_j)( 1 \, + \, o_{\l}(1)),
\end{equation}
where $o_{\l}(1) \to  0 \mbox{ if } \l \to + \infty$.
See please Lemma \ref{l:Ual} for the equation satisfied by $U_{A,\L}$. Furthermore we collected in the appendix some useful estimates on $\varphi_{a,\l}$. Such estimates are used in the expansion of the Euler Lagrange functional and its gradient in the neighborhood at infinity $V(m,\e)$. See  please \eqref{eq:Vme} for the definition of this set.

%%%%%%%%%%%%%%%%%%%%%%%%%%%%%
\section{ Lack of compactness and a neighborhood at infinity}
%%%%%%%%%%%%%%%%%%%%%%%%%%%%%

Our approach is variational. Therefore,  in order to detect critical points for the functional $J_{\varrho}$ with $\varrho = 8 \pi m$, we have to find all possible obstructions in deforming sublevel sets  $J_\varrho^a:=\{u:J_\varrho(u)\leq a\}$. To this aim we make use of a pseudogradient for $J_{\varrho}$,  constructed in Horak-Lucia \cite{Horak-Lucia} (see also \cite{Lucia, andrea}), whose flow lines  have  the following property:\\
For any  fixed initial data  $u_0$,  the flow line  $\eta(t,u_0)$ generated by this pseudogradient satisfies the following property:
\begin{itemize}
\item either $J_{\varrho}(\eta(t,u_0)) \to -\infty$ as $t \to + \infty$,
\item or  $\eta(t,u_0)$  enters any arbitrary small neighborhood of   the set of solutions $u_\beta$ of $(\mathcal{P}_{8m\pi - \beta})$ with $\beta \in [0,\ov{\e} )$ for some $\ov{\e} >0$.
\end{itemize}
Furthermore, for a solution $u_\beta$ with $ \beta >0$ we have the following alternative:  $(i)$ either $u_\beta$ converges to a solution $\ov{u}$ of $(\mathcal{P}_{8m\pi })$ as $\beta \to 0$   $(ii)$ or it has to blow up when $\beta \to 0$. In the latter  case,  it follows from  Proposition \ref{p:blowup} (see also  \cite{cl1, cl2}),  that the solutions $(u_\beta)$ have to belong to  some subset ${V}(m,\varepsilon)$, called in the sequel {\em neighborhood  at Infinity}, which is defined as:
\begin{eqnarray}\label{eq:Vme}
  {V}(m,\varepsilon) &:= \Big \{  u \in  \mathring{H}^1(\Sigma)\, : \, \|\nabla J_{\varrho} (u)\| < \varepsilon\, ; \, \, \exists \,  \lambda_1, \cdots, \lambda_m > {\varepsilon^{-1}} \mbox{ with } {\lambda_i} < C_1{\lambda_j}\, \,  \forall i \neq j; \nonumber\\
   &  \quad \exists \, a_1, \cdots, a_m \mbox{ with } d_g(a_i,a_j) \geq 2 \eta \, \, \forall i \neq j  \quad \mbox{such that } \| u -  \sum_{i = 1}^m  \varphi_{a_i,\lambda_i} \| < \varepsilon  \Big \},
\end{eqnarray}
where the space $  \mathring{H}^1(\Sigma)$ is defined in \eqref{eq:H10}, $\varepsilon $ is a small positive constant and  $C_1$, $\eta$ are fixed positive constants.\\
Hence, we are led to study the obstructions to decrease the functional $J_{\varrho}$ in the set $V(m, \varepsilon)$.
A first step consists in finding an appropriate parametrization of this set. To that aim,
following the ideas of A. Bahri and J.M. Coron   we consider the following minimization problem
\begin{equation}\label{f:min} \min_{ \alpha_i >0 ; a_i \in \Sigma; \lambda_i > 0}  \big\| u - \sum_{i=1}^m \alpha_i \varphi_{a_i,\lambda_i} \big\|\, .\end{equation}
We have the following Lemma whose proof is identical to the proof of Proposition 7 in  \cite{BC} (see also  Chen and Lin \cite[Lemma 3.2]{cl2}). Namely we have
\begin{lem}\label{minimi}
For $\varepsilon$ small, Problem~\eqref{f:min} has for any $u \in {V}(m,\varepsilon)$ only one solution (up to permutations on the indices). The variables $\alpha_i$'s satisfy $|\alpha_i -1| = O(\varepsilon) $.
\end{lem}
Hence every $u \in {V}(m,\varepsilon)$ can be written as
\begin{equation}\label{f:paramet}    u \, = \, \sum_{i=1}^m \alpha_i  \varphi_{a_i,\lambda_i} \, +  \, w,\end{equation}
where $\alpha_i$, $w$ satisfy
\begin{equation} \label{f:orthog} \begin{cases}
&  |\alpha_i -1| \leq c \varepsilon \quad \forall i,    \qquad    \|w \| \leq c \varepsilon , \\
&  \big <w,\varphi_{a_i,\lambda_i} \big >_g \,  =  \,  \big <w, {\partial \varphi_{a_i,\lambda_i}}/{\partial \lambda_i} \big >_g \, = \, 0 ,\, \quad\,  \big <w, { \partial \varphi_{a_i,\lambda_i}}/{\partial a_i} \big >_g \, = 0 \quad \forall\, \, i .
\end{cases}\end{equation}

In the following, for ${A}=(a_1,...,a_m)$ and $\Lambda=(\lambda_1,...,\lambda_m)$,
we denote
 \begin{equation} \label{eal}  E_{{A},\Lambda}^m:= \{ w\in \mathring{H}^1(\Sigma): w \mbox{ satisfies \eqref{f:orthog}} \}. \end{equation}
To keep the notation short we will write $\varphi_i$ instead of $\varphi_{a_i, \lambda_i}$.
The next Proposition  shows how to deal with the infinite dimensional part $w$ in the above representation:

\begin{pro}\label{p:wpart}\cite{ABL17}
Let $ u := \sum_{i =1}^m \alpha_i \varphi_i \in V(m,\varepsilon)$. Then there exists a unique $\overline{w} := \overline{w} (u) \in  E_{{A},\Lambda}^m$ such that:
$$  J_{\varrho} ( u + \overline{w}) \, = \,  \min \{ J_{\varrho} (u + w) \, : \,  w \in E_{{A},\Lambda}^m \}.$$
Furthermore, there exists a constant $C$ such that
\begin{equation} \label{eq:Estwbar} \| \overline{w} \| \, \leq \,  C  \sum_{i=1}^m   \Big(|\alpha_i-1|+ \frac{1}{\lambda_i}\Big) .\end{equation}
\end{pro}

%%%%%%%%%%%%%%%%%%%%%%%%%%%%%%%%%%%%%%%%
\section{ The expansion of the gradient in the neighborhood at infinity }
%%%%%%%%%%%%%%%%%%%%%%%%%%%%%%%%%%%%%%%%

In this section we provide an asymptotic expansion of the gradient of the Euler Lagrange functional $J_{\varrho}$ in the neighborhood at infinity $V(m,\e)$. For the sake of simplicity of the notation and since the variables $\l_i$'s are of the same order, we will write  in this section and in the sequel $O(1/\l^\g)$ instead of $\sum O(1/\l_k^\g)$.\\
In the first proposition we expand the gradient with respect to the concentration rates. Namely we have

\begin{pro}\label{derivelambda}
Let $u := \sum_{i=1}^m \a_i \varphi_{i} + w  \in V(m,\e)$ with $w\in E^m _{{A}, \Lambda}$ and $\varrho:= 8m\pi (1+\mu)$ with $\mu$ a small real number. It holds
\begin{equation*}
 \Big\langle \nabla J_{\varrho}(u), \lambda_i \frac{\partial  \varphi_i}{\partial  \lambda_i}  \Big\rangle_{g} = 16 \pi \alpha_i (\tau_i -\mu +\mu\tau_i)  - 64\pi^2 \sum_{j=1}^m  \frac{\ln\l_j}{\lambda_j^2}+  O\left( \sum |\alpha_k - 1|^2 + \|w\|^2 \right)+o\left(\frac{\ln\l}{\lambda^2} \right)
\end{equation*}
with
\begin{equation}\label{e:tau}   \tau_i \, := \, 1 \, - \, \frac{m \pi}{2\alpha_i-1} \,    \frac{\lambda_i^{4\alpha_i -2} \mathcal{F}^{A}_i(a_i) g^{A}_i(a_i)}{\int_{\Sigma} K e^u dV_g},
\end{equation}
where $\mathcal{F}^{A}_i$ and $g^{A}_i$  are defined in \eqref{eq:fkmi2}. Furthermore, we have the estimate
\begin{equation} \label{c:tau} | \tau_i | = O(\varepsilon + | \mu | ), \quad \forall i \in \{ 1, \cdots , m \}.\end{equation}
 \end{pro}

\begin{pf}
Recall that
\begin{equation}\label{g01}
 \Big\langle \nabla J_{\varrho}(u), h  \Big\rangle_{g} =  \Big\langle u , h    \Big\rangle_{g} - \frac{\varrho}{\int_\Sigma K e^u } \int_\Sigma K e^u h.
 \end{equation}
Using Lemma \ref{AA2} we get that
\begin{equation*}
 \Big\langle u , \lambda_i \frac{\partial  \varphi_i}{\partial  \lambda_i}  \Big\rangle_{g} = - 64\pi^2  \frac{\ln\l_i}{\lambda_i^2} \sum_{j=1}^m \a_j + 16 \pi \a_i
 + O( \frac{1}{\lambda ^2}) =  - 64\pi^2 m  \frac{\ln\l_i}{\lambda_i^2}  + 16 \pi \a_i  +  o( \frac{\ln\l}{\lambda ^2}) . \end{equation*}
For the other term of \eqref{g01}, using Lemma \ref{AA1}, we have
\begin{align*}
\int_\Sigma K e^u \lambda_i \frac{\partial  \varphi_i}{\partial  \lambda_i} dV_g & = \int_\Sigma K e^u \left( \frac{4}{1+ \l_i^2 \psi_i^2} - 8\pi  \frac{\ln\l_i}{\lambda_i^2}  + O( \frac{1}{\lambda ^2})\right)dV_g\\
& = -8\pi \frac{\ln\l_i}{\lambda_i^2} \int_\Sigma K e^u dV_g +  \int_{B_i} K e^u \frac{4}{1+ \l_i^2 \psi_i^2} dV_g+  O( \frac{1}{\lambda ^2})\int_\Sigma K e^u  dV_g .
\end{align*}
Now we need to estimate the second integral. Letting  $\ov{u}:= u-w$ and using Lemma A.4 in \cite{ABL17}, we have
\begin{align*}
\int_{B_i} K e^u \frac{4}{1+ \l_i^2 \psi_i^2}  & =  \int_{B_i} K e^{\ov{u}} \frac{4}{1+ \l_i^2 \psi_i^2}  +  \int_{B_i} K e^{\ov{u}}w \frac{4}{1+ \l_i^2 \psi_i^2}  +  \int_{B_i} K e^{\ov{u}}(e^w-1-w) \frac{4}{1+ \l_i^2 \psi_i^2}  \\
& = \int_{B_i} K e^{\ov{u}} \frac{4}{1+ \l_i^2 \psi_i^2} dV_g + O\bigg( (\sum \l_k^{4 \a_k - 2 })\| w \| \Big( \| w \| + \frac{1}{\l} + | \a_i-1|\Big)\bigg).
\end{align*}
Concerning the last integral, using Lemma \ref{eu}  we get
\begin{align*}
 \int_{B_i} K e^{\ov{u}} \frac{4}{1+ \l_i^2 \psi_i^2} dV_g & =   \Big(1+ 4\pi \sum_{j=1}^m  \frac{\ln\l_j}{\l_j^2}\Big) \int_{B_i} \frac{ 4 \lambda_i^{4\alpha_i}\mathcal{F}_i^{A}g_i^{A} e^{-u_{a_i}}}{(1+ \lambda_i^2| y_{a_i}(.) | ^2)^{2\alpha_i +1}} dV_{g_{a_i}}  + o\Big(\frac{\ln\l}{ \lambda^2}\Big) \int_{B_i} K e^{\ov{u}} \\
 & =  \Big(1+ 4\pi \sum_{j=1}^m  \frac{\ln\l_j}{\l_j^2}\Big) \frac{2\pi}{\a_i} \l_i^{4\a_i - 2}\mathcal{F}_i^{A}g_i^{A} (a_i) + o\Big(\frac{\ln\l}{ \lambda^2}\Big) \sum \l_k^{4\a_k - 2},
 \end{align*} by using
 $$  \int_0^{\l\eta} \frac{ 4\, r}{(1+r^2)^{2\a+1}} dr = \frac{1}{\a } + O\Big( \frac{1 }{\l^{4\a}}\Big)\quad \mbox{ and } \quad
 \int_0^{\l\eta} \frac{r^3}{(1+r^2)^{2\a+1}}dr = O \big(1\big).$$
 \noindent Thus we get
 \begin{align*}
\frac{\varrho}{\int_\Sigma K e^u } \int_\Sigma K e^u \lambda_i \frac{\partial  \varphi_i}{\partial  \lambda_i} dV_g = & \Big(1+ 4\pi \sum_{j=1}^m  \frac{\ln\l_j}{\l_j^2}\Big) 16\pi (1+\mu) \frac{2\a_i - 1 }{\a_i} (1- \tau_i)\\
& +  o\Big(\frac{\ln\l}{ \lambda^2}\Big) + O\Big( \| w \| ^2 + | \a_i - 1 | ^2 \Big).\end{align*}
 The result follows by summing the previous estimates.
\end{pf}

In the next proposition we expand the gradient $\n J_{\varrho}$ with respect to the gluing parameters $\a_i$'s.
\begin{pro}\label{derivealpha}
Let $u := \sum_{i=1}^m \a_i \varphi_{i} + w  \in V(m,\e)$ with $w\in E^m _{{A}, \Lambda}$ and $\varrho:= 8m\pi (1+\mu)$ with $\mu$ a small real number. It holds
\begin{align*}
 \Big\langle \nabla J_{\varrho}(u),  \varphi_i \Big\rangle_{g} = & 32 \pi \ln \l_i \big( (\a_i-1) + (\tau_i - \mu +\mu \tau_i) \big)  +  O\left( |\alpha_i - 1| + \|w\|^2\ln\l \right)\\
 &  +  O\left( \| w \| (1+|\a_i-1| \ln\l ) + \sum | \tau_k -\mu +\mu \tau_k| +\frac{\ln^2\l}{\lambda^2} + \frac{\ln\l_i}{\lambda_i^{4\a_i-2}}\right).
\end{align*}
\end{pro}
\begin{pf}
Using Lemma \ref{AA2} we get that
 \begin{align*}
 \Big\langle u ,  \varphi_i \Big\rangle_{g} & = \a_i \Big( 32 \pi   \ln\l_i  + 64 \pi^2 H(a_i,a_i) - 16 \pi \Big) + 64 \pi^2 \sum_{j\neq i} \a_j G(a_i,a_j)  + O\Big( \frac{\ln\l}{\l^2} \Big) ,\\
\int_\Sigma K e^u  \varphi_i dV_g & = \int_\Sigma K e^{\ov{u}}  \varphi_i dV_g + \int_\Sigma K e^{\ov{u}}  w\varphi_i dV_g + \int_\Sigma K e^{\ov{u}} (e^w - 1 - w ) \varphi_i dV_g\\
& =  \int_\Sigma K e^{\ov{u}}  \varphi_i dV_g  + O\Big( \| w \| \Big(1 + | \a_i - 1 | \ln \l_i  + \| w \| \ln\l_i \Big) \sum \l_k^{4\a_k -2} \Big) .
\end{align*}
{
For the last integral we notice that it follows from Lemma \ref{AA1} that  $\varphi_i$ and $e^{\ov{u}} $ are bounded outside the balls $B_k$. Moreover inside  each  ball $B_k$, for $k\neq i$, using Lemmas \ref{AA1} and \ref{eu}, we derive
}
\begin{align*}
\int_{B_k} K e^{\ov{u}}  \varphi_i dV_g & = 8\pi \int_{B_k} \frac{  \lambda_k^{4\alpha_k}\mathcal{F}_k^{A}g_k^{A} e^{-u_{a_k}}}{(1+ \lambda_k^2| y_{a_k}(.) | ^2)^{2\alpha_k}} G(a_i,.) dV_{g_{a_k}} + O\Big( \frac{\ln\l }{\l^2}\Big)\int K e^ {\ov{u}}\\
& = 8\pi   \lambda_k^{4\alpha_k-2}\mathcal{F}_k^{A}g_k^{A}(a_k)G(a_i,a_k) \frac{\pi}{2\a_k -1} + O\Big(\Big( \frac{\ln\l }{\l^2} + \sum \frac{ | \a_j -1 | }{\l^{3/2}} \Big)\sum  \lambda_j^{4\alpha_j-2}\Big)\end{align*}
by using the fact that (using \eqref{e:nxi})
\begin{align*}
\int_0^{\l\eta} \frac{r^3 dr}{(1+r^2)^{2\a} }  = \int_0^{\l\eta} \frac{r^3 \xi(r) dr}{(1+r^2)^{2} } & = \int_0^{\l\eta} \frac{r^3dr}{(1+r^2)^{2} }  + O\Big( |\a-1| \int_0^{\l\eta} \frac{r^3\sqrt{r}dr}{(1+r^2)^{2} }\Big) \\
& =  O\Big( \ln \l +  |\a-1|\sqrt{\l}\Big).  \end{align*}
In the ball $B_i$, it holds
\begin{align*}
\int_{B_i}  K & e^{\ov{u}}  \varphi_i dV_g \\
& = \int_{B_i} \frac{  \lambda_i^{4\alpha_i}\mathcal{F}_i^{A}g_i^{A} e^{-u_{a_i}}}{(1+ \lambda_i^2| y_{a_i} | ^2)^{2\alpha_i}}\Big( 4\ln\l_i - 2 \ln(1+ \lambda_i^2| y_{a_i} |^2) + 8\pi H(a_i,.)\Big) dV_{g_{a_i}}  + O\Big( \frac{\ln\l }{\l^2}\Big)\int K e ^{\ov{u}}\\
& = \Big( 4 \ln\l_i + 8\pi H(a_i,a_i) \Big) \Big( \frac{\pi}{2\a_i -1} + O\Big( \frac{1}{\l_i^{4\a_i -2}} + \frac{\ln\l}{\l^2} + \frac{ | \a_i -1|}{\l^{3/2}}\Big)\Big)  \lambda_i^{4\alpha_i-2}\mathcal{F}_i^{A}g_i^{A} (a_i)\\
& - 2 \int_{B_i} \frac{  \lambda_i^{4\alpha_i}\mathcal{F}_i^{A}g_i^{A} e^{-u_{a_i}}}{(1+ \lambda_i^2| y_{a_i} | ^2)^{2\alpha_i}} \ln(1+ \lambda_i^2| y_{a_i} |^2) dV_{g_{a_i}}  + O\Big( \frac{\ln\l }{\l^2}\Big)\int K e ^{\ov{u}}.
\end{align*}
Observe that
\begin{align*}
& \int_0^{\l\eta} \frac{2r}{(1+r^2)^{2\a}} \ln(1+r^2) dr = \frac{1}{(2\a -1)^2} + O\Big( \frac{\ln \l }{\l^{4\a-2}}\Big),\\
& \int_0^{\l\eta} \frac{r^3}{(1+r^2)^{2\a}} \ln(1+r^2) dr \leq c \, \ln\l  \int_0^{\l\eta} \frac{r^3}{(1+r^2)^{2\a}} dr \leq c \, \ln \l \Big( \ln \l +  |\a-1|\sqrt{\l}\Big).
\end{align*}
Thus we obtain
\begin{align*}
\int_{B_i} K e^{\ov{u}}  \varphi_i dV_g = & \Big( 4 \ln\l_i + 8\pi H(a_i,a_i) \Big)  \frac{\pi}{2\a_i -1} \lambda_i^{4\alpha_i-2}\mathcal{F}_i^{A}g_i^{A} (a_i) - \frac{2\pi}{(2\a_i - 1)^2}\lambda_i^{4\alpha_i-2}\mathcal{F}_i^{A}g_i^{A} (a_i)\\
&  + O\Big(\ln\l +  \Big( \frac{\ln^2\l}{\l^2}+|\a_i-1|\frac{\ln\l}{\l^{3/2}}
\Big)\sum \lambda_j^{4\alpha_j-2}\Big).\end{align*}
As in the proof of Proposition \ref{derivelambda}, summing the previous estimates, we derive the result.
\end{pf}

\noindent
Combining Propositions  \ref{derivelambda} and \ref{derivealpha}, we obtain
\begin{cor}\label{c:al}
Let $u := \sum_{i=1}^m \a_i \varphi_{i} + w  \in V(m,\e)$ with $w\in E^m _{{A}, \Lambda}$ and $\varrho:= 8m\pi (1+\mu)$ with $\mu$ a small real number. It holds
\begin{align*}
 \Big\langle & \nabla J_{\varrho}(u),  \frac{\varphi_i}{\ln\l_i} - \frac{2}{\a_i} \lambda_i \frac{\partial  \varphi_i}{\partial  \lambda_i}  \Big\rangle_{g} \\
 & =   32 \pi  (\a_i-1)  +  O\left( \|w\|^2 + \frac{ \| w \| }{\ln\l } +  \frac{ 1 }{\ln\l } \sum (| \tau_k -\mu +\mu \tau_k| +| \a_k - 1 |) +\frac{\ln\l}{\lambda^2} + \frac{1}{\lambda^{4\a_i-2}}\right).
\end{align*}
\end{cor}

\begin{lem}\label{sumtau} Let $u := \sum_{i=1}^m \a_i \varphi_{i} + w  \in V(m,\e)$ with $w\in E^m _{{A}, \Lambda}$ and $\varrho:= 8m\pi (1+\mu)$ with $\mu$ a small real number. Let $\tau_i$ be defined in \eqref{e:tau} and assume that $\sum  |\a_i -1| \ln \l_i  $ is small. Then, it holds
 $$ \sum_{i=1}^m \tau_i = \frac{\pi}{2}\frac{m \ln\l_1}{\int_{\Sigma} K e^u dV_g} \mathcal{L}(A)  + 4\pi m \sum_{j=1}^m \frac{\ln\l_j }{\l_j^2} + o\Big(  \frac{\ln\l }{\l^2} \Big) + \sum O\Big( | \a_k - 1 |^2\Big).$$

\end{lem}
\begin{pf}
In this case, since we assumed that $\sum |\a_i -1| \ln \l_i $ is small, we derive that $\l_i^{4(\a_i -1)} = 1 +O(| \a_i -1 | \ln \l_i)$ for each $i$. Moreover we have
\begin{equation}\label{e:w1} \int_\Sigma K e^{\ov{u}} w dV_g + \int_\Sigma K e^{\ov{u}} (e^w -1-w) dV_g = O\Big(\Big( \| w \|^2 + \sum | \a_k -1 |^2 + \frac{1}{\l^2}\Big)\sum\l_k^{4\a_k -2} \Big).\end{equation}
Hence the result follows from the previous estimate and  \eqref {eq:Raffaella2}.
\end{pf}

Lastly arguing as above, we derive the following asymptotic expansion of the gradient with respect to the concentration points $(a_1, \cdots,a_m)$. Namely we prove that
\begin{pro}\label{deriva}
Let $u := \sum_{i=1}^m \a_i \varphi_{i} + w  \in V(m,\e)$ with $w\in E^m _{{A}, \Lambda}$ and $\varrho:= 8m\pi (1+\mu)$ with $\mu$ a small real number. It holds
\begin{align*}
 \Big\langle \nabla J_{\varrho}(u), & \frac{1}{\l_i} \frac{\partial \varphi_i}{\partial a_i} \Big\rangle_{g} \, = \, -8 \pi(1 + \mu) \frac{\n  \mathcal{F}^A_i(a_i)}{\l_i}   +  O\left(   \frac{| \mu | }{\l} + \frac{|\tau_i|}{\l } + \frac{\ln\l }{\l ^2}+ \sum_{k=1}^{m} |\alpha_k - 1|^2 + \|w\|^2 \right) .
\end{align*}
\end{pro}

%%%%%%%%%%%%%%%%%%%%%%%%%%%%%%%%%%%%%%%%%

\section{ Critical points at infinity and their indices }

Critical points at infinity of the functional $J_{\varrho}$ are accumulation points of some orbits of the negative gradient flow $- \n J_{\varrho}$ which enter some $V(m,\e)$ and remain there indefinitely. In this section we recall for the sake of completeness the full characterization of these critical points proved in \cite{ABL17} and restate their contribution to the difference of topology between the level sets of the functional $J_{\varrho}$.\\
We start with an expansion of $J_{\varrho} $ in the neighborhood at infinity $V(m,\e)$.

%%%%%%%%%%%%%%%%%%%%%%%%%%%%%%%%%%%%%%%%%

\begin{pro}\label{devJ}
Let $u := \sum_{i=1}^m \a_i \varphi_{i} + w \in V(m,\e)$ with $w\in E^m _{{A},\Lambda}$. Assume that $| \a_i -1 | \ln\l_i$ is small for each $i$ and \eqref{eq:Estwbar} holds. Then
\begin{align*}
J_{8m\pi} & (u)  =  -8 \pi m (1+\ln(m\pi))-8\pi  \mathcal{F}_m^K(a_1,...,a_m) -4\pi \sum_{i=1}^m | \tau'_i |^2  + 16\pi \sum_{i=1}^m (\alpha_i-1)^2  \ln \lambda_i \\
 & -  4\pi \frac{\ln\l_1}{\l_1^2 \mathcal{F}^{A} _1 (a_1)}  \sum_{i=1}^m \Big( \Delta \mathcal{F}^{A} _i (a_i)  -2 K_{g}(a_i)\mathcal{F}^{A} _i (a_i)\Big) +  \sum_{k=1}^m  \Big\{ O \Big( (\alpha_k -1)^2   +\frac{1}{\lambda_k^2}   \Big) + o(| {\tau}'_k |^2 ) \Big\} ,
\end{align*}
 where $\mathcal{F}^{A}_i$ and $g^{A}_i$  are defined in \eqref{eq:fkmi2} and
$$ {\tau}'_i = 1- \frac{ {m}\,  \lambda_i^{4\alpha_i-2} \mathcal{F}^{A}_i(a_i)}{\sum_{k=1}^m \lambda_k^{4\alpha_k-2} \mathcal{F}^{A}_k (a_k)}.$$
\end{pro}

\begin{pf}   The proof follows from Lemmas \ref{AA2} and \ref{intkeu} and the following computations. Let us denote by
$$\Gamma := {\sum_{i=1}^m \frac{\lambda_i^{4 \alpha_i-2} \mathcal{F}^{A}_i (a_i)g_i^{A}(a_i)}{2 \alpha_i-1}} \qquad \mbox{ and } \qquad  \tilde{\tau}_i = 1- \frac{ \frac{m}{2\alpha_i-1}\lambda_i^{4\alpha_i-2} \mathcal{F}^{A}_i(a_i)g_i^{A}(a_i)}{\sum \frac{1}{2\alpha_k-1}\lambda_k^{4\alpha_k-2} \mathcal{F}^{A}_k (a_k)g_k^{A}(a_k)}.$$
Since $w\in E^m _{{A},\Lambda}$, then \eqref{e:w1} holds and using Lemma \ref{intkeu}, we have
\begin{align*}
  \ln \Big(\int_{\Sigma} K e^u  \Big)  & =  \ln(\pi \Gamma) +  \ln\Big\{ 1+  \frac{1}{2\Gamma} \sum_{i=1}^m \Big( \Delta \mathcal{F}^{A} _i (a_i)  -2 K_{g}(a_i)\mathcal{F}^{A} _i (a_i)\Big)\ln \lambda_i \\
  & + \frac{4\pi}{\Gamma} \Big( \sum_{j=1}^m\frac{\ln\l_j}{\l_j^2} \Big) \sum_{i=1}^m   \lambda_i^{4\alpha_i-2} \mathcal{F}_i^{A}(a_i)+ O \big( \sum_{k=1}^m  \big \{ |\alpha_k - 1|^2 + \frac{1}{\lambda^2_k}  \big \} \big) \Big\} \\
 &= \ln \pi  -\frac{1}{m} \sum_{i=1}^m \ln(1-\tilde\tau_i)  + \frac{1}{2\Gamma} \sum_{i=1}^m \Big( \Delta \mathcal{F}^{A} _i (a_i)  -2 K_{g}(a_i)\mathcal{F}^{A} _i (a_i)\Big)\ln \lambda_i  \\
 & + \frac{1}{m} \sum_{i=1}^m \ln \Big(\frac{m \lambda_i^{4\alpha_i-2} \mathcal{F}^{A}_i(a_i) g_i^A(a_i)}{2\alpha_i-1}  \Big) + \frac{4\pi}{\Gamma} \Big( \sum_{j=1}^m\frac{\ln\l_j}{\l_j^2} \Big) \sum_{i=1}^m   \lambda_i^{4\alpha_i-2} \mathcal{F}_i^{A}(a_i)\\
&  + O \Big( \sum_{k=1}^m  \big \{ |\alpha_k - 1|^2 + \frac{1}{\lambda^2_k}  \big \} \Big)  \nonumber\\
  & =  \ln \pi +  \frac{1}{m} \sum_{i=1}^m \Big\{ \tilde\tau_i + \frac{\tilde\tau_i^2}{2}   \Big\} + \frac{1}{2\Gamma}  \sum_{i=1}^m  \Big( \Delta \mathcal{F}^{A} _i (a_i)  -2 K_{g}(a_i)\mathcal{F}^{A} _i (a_i)\Big)\ln \lambda_i  \\
& + \ln m + \frac{1}{m} \sum_{i=1}^m \{- \ln(2\alpha_i-1)   +   (4\alpha_i - 2)\ln \lambda_i +  \ln(\mathcal{F}^{A}_i(a_i)) + \ln(g_i^A(a_i)) \}\nonumber \\
 & + \frac{4\pi}{\Gamma} \Big( \sum_{j=1} ^m\frac{\ln\l_j}{\l_j^2} \Big) \sum_{i=1}^m   \lambda_i^{4\alpha_i-2} \mathcal{F}_i^{A}(a_i)
   +  O  \Big(\sum_{k=1}^m  \Big\{ |\alpha_k -1|^2 +|\tilde\tau_k|^3 + \frac{1}{\lambda_k^2} \Big\}  \Big).
\end{align*}

We remark that $\sum_{i=1}^m \tilde\tau_i =0$ and for each $i$, we have $ \tilde\tau_i  = O(\e)$ which implies that $\l_i^2 \mathcal{F}^{A} _i (a_i) = \l_j^2 \mathcal{F}^{A} _j (a_j)(1+ O(\sum | \a_k - 1 | \ln\l_k))$ for each $i,j$ and $\ln\l_i = \ln\l_1 + O(1)$. Furthermore, we have
 \begin{align*}
 &  \ln \big( 2\alpha_i-1 \big) = \ln \big( 1+2(\alpha_i-1) \big)=2(\alpha_i-1) + O(|\alpha_i-1|^2),\\
 & g_i^A(a_i)= 1+ O \Big(\sum_{k=1}^m |\alpha_k -1| \Big) \quad ; \quad   \Gamma = \sum_{i=1}^m \lambda_i^2 \mathcal{F}^{A}_i(a_i)  + O \Big(\sum_{k=1}^m |\alpha_k-1| \lambda_k^2 \ln \lambda_k \Big)\\
 & \mbox{and} \quad | \tilde{\tau}_i |^2  = | \tau'_i |^2  + o( | \tau'_i |^2 ) + \sum  O(| \a_k - 1 |^2) .\end{align*}
 Hence, the result follows.
\end{pf}

 In \cite{ABL17} we constructed the following  decreasing pseudogradient of the Euler Lagrange functional $J_{\varrho}$ in the neighborhood at Infinity $V(m,\e)$:
\begin{pro}\cite{ABL17}\label{p:champ}
Let $\varrho = 8 \pi m$ with $m\geq 1$ and assume that
the function $K$ satisfies the condition $(i)$ of $(\mathcal{N}_m)$. Then there exists a pseudogradient
$W$ defined in $V(m,\e)$ and satisfying the following properties:
There exists a constant $C$ independent of $ u = \sum_{i=1}^m \a_i \varphi_{a_i,\l_i}+\bar{w}$ such that
\begin{enumerate}
\item[{\rm (1)}]
$\displaystyle{ \langle-\n J_{\varrho}(u),W\rangle \, \geq \,
 C \sum_{i=1}^m \Big( |\a_i - 1| \, + \, |\tau_i|  \, + \, \frac{|\n \mathcal{F}^{A}_i(a_i)|}{\l_i}  \, + \,  \frac{\ln\l_i}{\l_i^2} \Big) }$,
\item[{\rm (2)}] $\displaystyle{ \langle-\n J_{\varrho}(u), W + \frac{\partial \ov{w}(W)}{\partial(\a,\l,a)}\rangle \, \geq \, C \sum_{i=1}^m
\Big( |\a_i - 1| \, + \, |\tau_i|  \, + \, \frac{|\n \mathcal{F}^{A}_i(a_i)|}{\l_i}  \, + \,  \frac{\ln\l_i}{\l_i^2} \Big)  }$,
\item[{\rm (3)}]
$|W|$ is bounded and the only region where the variables $\l_i$'s increase along the flow lines of $W$ is the region where
$(a_1,\cdots,a_m)$ is very close to  a critical point $q:=(q_1,\cdots,q_m)$ of $\mathcal{F}^K_m$ with $ q \in\mathcal{K}^-_m$.
\end{enumerate}
\end{pro}

Following the program developed in \cite{ABL17}, we deduce  from  Proposition \ref{p:champ} the following characterization of the critical points at infinity:

\begin{pro}\cite{ABL17}\label{c:cpatinfinity} Let $\varrho= 8 \pi m$, $m\geq 1$.
The critical points at Infinity of $J_{\varrho}$  are in  one to one correspondence with critical points $Q:=(q_1,\cdots,q_m)$ of $\mathcal{F}_m^{K}$  satisfying $Q \in \mathcal{K}^-_m, $ that   will be denoted   by $(Q)_{\infty}$.
Furthermore  the energy level of such a critical point at Infinity $(Q)_{\infty}$ denoted $C_{\infty}(Q)_{\infty}$ is given by:
$$ C_{\infty}( Q )_{\infty} \, =  \, - 8 \pi m (1+\ln(m\pi)) \, - 8\pi \mathcal{F}^{K}_m ( Q) . $$
Moreover the Morse index of such a critical point at Infinity $(Q)_{\infty} $ is given by:
$$ \i _\infty (Q)  := \, 3m - 1 - Morse(\mathcal{F}^{K}_m,Q).$$
\end{pro}

\bigskip
\noindent
We point out that  around \emph{a critical point at infinity} there is  a Morse type reduction. Indeed denoting by
\begin{align}\label{eq:vmq}
V(m,Q,\e) \,:  = \{ u= \sum_{i=1}^{m} \a_i \varphi_{a_i,\l_i} \, + w \,  \in V(m,\e): \,  \forall i =1,\cdots,m, \, |a_i-q_i| < \e; \, \nonumber  \\
 |\a_i - 1| \ln\l_i <\e,  \mbox{ and }
\forall i \neq j \quad  | \frac{\l_i^2  \mathcal{F}^A_i(a_i)}{ \l_j^2  \mathcal{F}^A_j(a_j)} \, -  1 | < \e \},
\end{align}
where $Q:=  (q_1, \cdots,q_m)$ is a critical point of $\mathcal{F}^K_m$ with $\mathcal{L}(Q) < 0$, we have the following  { \it  Morse Lemma at Infinity},whose proof is inspired by the proof of a similar statement for Yamabe type flows written down by A. Bahri  in \cite{Binvariant} (pages 415-417). We state our result as follows:

\begin{lem}\label{eq:morselem}
Let $ u= \sum_{i=1}^{m} \a_i \varphi_{a_i,\l_i} \, + w \,  \in V(m,Q,\e)$, then there exists a change of variables
$$ ( a_1,\cdots,a_m, \alpha_1,\cdots,\alpha_m,  \lambda_1,\cdots,\lambda_m,w) \mapsto  ( \ov{a}_1,\cdots, \ov{a}_m, \b_1,\cdots,\b_m,\L_1, y_2,\cdots,y_m,V) $$ (where $(\ov{a}_1,\cdots, \ov{a}_m)$ is close to $Q$, the $\b_i$'s and the $y_i$'s are small and $\L_1$ is very large) such that
$$  J_{\varrho}(u)   =   -8 \pi m (1+\ln(m\pi))-8\pi \mathcal{F}_m^K(\ov{a}_1,...,\ov{a}_m) + \sum_{i=1}^m \b_i^2  - \sum_{i=2}^m y_i^2 +  (4\pi -\s) \frac{ (-\mathcal{L} (Q)) }{ \mathcal{F}^{Q} _1 (q_1)} \frac{\ln {\L}_1}{\L_1 ^2} + \| V \|^2$$
where $\s$ is a small positive constant.
\end{lem}

\begin{pf}  We recall that  for $ \sum \a_i \varphi_i \in V(m,\e)$, we have by  Proposition \ref{p:wpart} that  there exists a unique  $\ov{w}$ which minimizes $J_\varrho(\sum \a_i \varphi_i + w)$ in the space $E_{ A,\Lambda}^m$. Hence, by the classical Morse Lemma, there exists a change of variable $w -\ov{w} \to V$ so that
\begin{equation}\label{ff11} J_\varrho(u)  = J_\varrho(\ov{u}) + \| V \|^2  \quad\mbox{ where }  {u} := \sum \a_i \varphi_i  + {w} \quad \mbox{ and } \quad \ov{u} := \sum \a_i \varphi_i  + \ov{w} .\end{equation}
Furthermore
for $\e' >0$ small (with  $\e / \e'$ small) and   $W$    the pseudogradient   constructed in Proposition \ref{p:champ}, we have that    for $\ov{u} := \sum \a_i \varphi_i + \ov{w}  \in V(m, \e')$ there holds:
\begin{equation}\label{**1*} \langle \n J_\varrho(\ov{u}) , W \rangle \leq - c \sum \Big( | \a _i - 1 | + | \tau_i | + \frac{| \n \mathcal{F}_i^A (a_i) | }{\l_i} + \frac{\ln\l_i}{\l_i ^2} \Big) .\end{equation}
Next  let $\s >0$ be a small constant and set
\begin{align}  I(\ov{u} )   := & -8 \pi m (1+\ln(m\pi))-8\pi  \mathcal{F}_m^K(a_1,...,a_m) - (4\pi -\s)  \sum_{i=1}^m | {\tau}'_i|^2   \label{ff12} \\
&  + (16\pi +\s) \sum_{i=1}^m (\alpha_i-1)^2  \ln \lambda_i ^{2\a_i -1}  -  (4\pi -\s) \frac{\ln\l_1^{2\a_1-1}}{\l_1^{4\a_1 -2} \mathcal{F}^{Q} _1 (q_1)} \mathcal{L} (Q) . \nonumber
\end{align}
where $Q:= (q_1,\cdots,q_m)$. Since we assumed that $| \a_i -1| \ln\l_i$ is small for each $i$, it is easy to see that
\begin{equation} \label{**3*}0 <  I(\ov{u}) - J_\varrho(\ov{u})  \leq  2\s \sum \Big(  (\alpha_i-1)^2  \ln \lambda_i  +  \frac{\ln\l_1}{\l_1^2 } +  | {\tau}'_i |^2 \Big) \quad \mbox{ for each } \ov{u}.\end{equation}
Furthermore it follows from Proposition \ref{p:champ} that
\begin{equation}\label{**2*} \langle \n I (\ov{u}) , W \rangle \leq - c \sum \Big( | \a _i - 1 | + | \tau_i | + \frac{| \n \mathcal{F}_i^A (a_i)| }{\l_i} + \frac{\ln\l_i}{\l_i ^2} \Big) .\end{equation}
Next for  $\ov{u}_0 \in V(m,\e)$ we consider  the following  differential equation
\begin{equation}\label{flot} \frac{\partial u}{\partial s} = W(u) \quad ; \quad u(0)= \ov{u}_0,\end{equation}
whose  solution is $h_s(\ov{u}_0)$ where $h_s$ is the 1-parameter group generated by $W$.

Note that, for $\ov{u} := \sum \a_i \varphi_{a_i,\l_i} + \ov{w}$  as far as $ h_s(\ov{u}) \in V(m,Q,\e)$ we have  that
$$ h_s(\ov{u}) = \sum \a_i(s) \varphi_{a_i(s), \l_i(s)} + \ov{w(s)},$$
that is $ \ov{w(s)}$ satisfies conclusions of  Proposition \ref{p:wpart}.\\
Next we  claim:\\
{\bf CLAIM 1:} There exists $\ov{s} > 0$ such that $I(h_{\ov{s}} (\ov{u}_0) ) = J_\varrho(\ov{u}_0)$.\\
Observe that $I(h_s(\ov{u}_0))$ is a decreasing function with respect to $s$. Hence there exists at most  a unique solution to the equation $I(h_s(\ov{u}_0)) = J_\varrho(\ov{u}_0)$.
The only cases where there could be no solution are
\begin{itemize}
\item either $h_s(\ov{u}_0)$ exits $V(m,\e')$ (outside this set, we loose \eqref{**1*} since $W$ is defined only in $V(m,\e')$) before reaching this level.
\item  or $h_s(\ov{u}_0)$ will build a critical point at infinity before reaching the level $J_\varrho(\ov{u}_0)$.
\end{itemize}
We will prove that these two cases cannot occur. In fact, for the first one, since $\ov{u}_0 \in V(m,\e)$ then, to exit  $V(m,\e')$,  the flow line has to travel from $V(m, \e'/2)$ to the boundary of $V(m,\e')$. Note that, using \eqref{**2*},
we have that: $\partial I(h_{{s}} (\ov{u}_0) ) / \partial s \leq -c(\e')$ along this path, independent of $\e$, but depending on $\e'$. Also, by \eqref{**2*}, the time to travel from  $V(m, \e'/2)$ to the boundary of $V(m,\e')$ is lower-bounded by a constant $c'(\e')$, because $| W|$ is bounded and the distance to travel is lower-bounded by a constant $c>0$. Therefore, $I(h_{{s}} (\ov{u}_0) )$ decreases at least by $c(\e') c'(\e')$ during this trip.\\
However, using \eqref{**3*}, since $\ov{u}_0 \in V(m,\e)$, it follows that  $ J_\varrho(\ov{u}_0) < I(\ov{u}_0) \leq J_\varrho(\ov{u}_0) + c(\e)$. Hence we have  to choose $\e$ small with respect to $\e'$ so that $I(h_{{s}} (\ov{u}_0) )$ reaches the level $J_\varrho(\ov{u}_0) $ before leaving the set $V(m, \e')$.\\
Concerning the second case, the flow line $h_{{s}} (\ov{u}_0) $ will enter $V(m,\e_1)$ for each $\e_1>0$. Observe that
$$  J_\varrho (h_s (\ov{u}_0 )) < J_\varrho (\ov{u}_0 ) \quad ; \quad  0 < I(h_{{s}} (\ov{u}_0) ) - J_\varrho(h_{{s}} (\ov{u}_0) ) \to 0 \quad \mbox{ as } s \to \infty \, \, (\mbox{for }  \e_1 \to 0) .$$
Hence $I(h_{{s}} (\ov{u}_0) )$ has to reach the level $J_\varrho (\ov{u}_0 ) $ and therefore this case cannot occur. \\
Our claim is thereby proved.\\
Conversely, taking $\e'' > 0$ small with respect to $\e$ and given $\ov{u}'_0 \in V(m,\e'')$, arguing by the same way (and using $-W$ as an increasing pseudogradient) : there exists $\ov{s}' > 0$ such that $J_\varrho(h_{-\ov{s}'}( \ov{u}'_0)) = I(\ov{u}'_0)$. \\
Hence, $h_s$ is the required isomorphism. \\
It follows from \eqref{ff12} and the previous claim that
\begin{align*}
J_{\varrho}(\sum_{i=1}^m\alpha_i\varphi_{a_i, \l_i}\, + \, \ov{w} )   := & -8 \pi m (1+\ln(m\pi))-8\pi \mathcal{F}_m^K(\ov{a}_1,...,\ov{a}_m)   - (4\pi -\s)  \sum_{i=1}^m | {\ov{\tau}}'_i |^2 \\
& + (16\pi +\s) \sum_{i=1}^m (\ov{\alpha}_i-1)^2  \ln (\ov{\lambda}_i) ^{2\ov{\a}_i-1}  -  (4\pi -\s) \frac{\ln(\ov{\l}_1)^{2\ov{\a}_1-1}}{(\ov{\l}_1)^{4\ov{\a}_1-2} \mathcal{F}^{Q} _1 (q_1)} \mathcal{L} (Q) ,
\end{align*}
where $\ov{\tau}'_i$ has the same definition as $\tau'_i$ (see Proposition \ref{devJ}) using the new variables.\\
In order to achieve the split of variables as claimed in the Lemma we need to perform some changes of variables.
We first consider the  following  change of variables
\begin{align*} \psi_1:  Y:= &  (\ov{a}_1,\cdots,\ov{a}_m,\ov{\a}_1,\cdots,\ov{\a}_m,\ov{\l}_1,\cdots,\ov{\l}_m) \mapsto (\ov{a}_1,\cdots,\ov{a}_m,\ov{\a}_1,\cdots,\ov{\a}_m ,\L_1,\cdots,\L_m)\\
& \mbox{ with }\L _i := \ov{\l}_i^{2\a_i-1} \end{align*}
By computing the determinant $ \det  (D \psi_1(Y)) $, it is easy  to show that $ \psi _1$ is a diffeomorphism. With this change of variables, the functional reads as follows:
\begin{align*}
J_{\varrho}(\sum_{i=1}^m\alpha_i\varphi_{a_i, \l_i}\, + \, \ov{w} )   := & -8 \pi m (1+\ln(m\pi))-8\pi \mathcal{F}_m^K(\ov{a}_1,...,\ov{a}_m)   - (4\pi -\s)  \sum_{i=1}^m | {{\tau}}''_i |^2 \\
& + (16\pi +\s) \sum_{i=1}^m (\ov{\alpha}_i-1)^2  \ln {\Lambda}_i  -  (4\pi -\s) \frac{\ln {\L}_1}{{\L}_1^2 \mathcal{F}^{Q} _1 (q_1)} \mathcal{L} (Q) ,
\end{align*}
where $$ {\tau}''_i = 1- \frac{ {m}\,  \Lambda_i^{2} \mathcal{F}^{\ov{A}}_i(\ov{a}_i)}{\sum_{k=1}^m \Lambda_k^{2} \mathcal{F}^{\ov{A}}_k (\ov{a}_k)}.$$
Furthermore we perform  the following change of variables:
$$ \psi_2:  Y:= (\ov{a}_1,\cdots,\ov{a}_m,\ov{\a}_1,\cdots,\ov{\a}_m,{\L}_1,\cdots,{\L}_m) \mapsto (\ov{a}_1,\cdots,\ov{a}_m,{\b}_1,\cdots,{\b}_m ,x_1,\cdots,x_m)$$
{ with }$$ \b _i := \sqrt{16\pi +\s} (\ov{\a}_i -1) \sqrt{\ln\L_i}\, \, \forall\, \,  1\leq i \leq m \quad ; \quad x_1:= {\L}_1 \quad ; \quad x_i := \sqrt{4\pi -\s}  \tau''_i \, \, \forall \, \, 2\leq i \leq m. $$

%%%%%%%%%%%%%%%%%%%%%%%

Next we  claim:\\
{\bf CLAIM 2:}  $\psi_2$ is a diffeomorphism on the set $(\a, a, \l)$ satisfying the conditions in the definition of $V(m,\e,Q)$, see \eqref{eq:vmq}.\\
Indeed  let $H:= (b_1,\cdots,b_m, \g_1,\cdots,\g_m,\xi_1,\cdots,\xi_m)$. We need to solve $D\psi_2(Y)(H)=0$ and to prove that the unique solution is $H=0$. Let $\psi_2^j$ be the $j$-th component of $\psi_2$. It holds
 $$ D\psi_2(Y)(H)=0 \Longleftrightarrow D \psi_2^j(H) = 0 \, \,  \forall \, \, 1 \leq j \leq 3m \Longleftrightarrow \begin{cases} b_j = 0 \, \,  \forall \, \, j=1,\cdots, m,\\
 \sqrt{\ln\L_i} \g_i + \frac{(\ov{\a}_i - 1) }{2 \L_i \sqrt{\ln \L_i}} \xi_i = 0 \, \,  \forall \, \, i \leq  m,\\
 \xi_1 = 0\\
 \sum_{j=2}^m \frac{\partial x_i}{\partial \L_j} \xi_j = 0 \, \,  \forall \, \, i=2,\cdots, m. \end{cases}$$
 We will focus on the last equation. Let $\mathcal{F}_k:= \mathcal{F}_k^{\ov{A}}(\ov{a}_k)$ and  $D:= \sum_{k=1}^m \L_k^2 \mathcal{F}_k(a_k)$.  Observe that
 $$ \L_i \frac{\partial x_i}{\partial \L_i} = -\frac{2m}{D^2} \L_i ^2 \mathcal{F}_i \Big( \sum_{k \neq i} \L_k ^2 \mathcal{F}_k\Big) \quad ; \quad  \L_j \frac{\partial x_i}{\partial \L_j} = \frac{2m}{D^2} \L_i ^2 \mathcal{F}_i \L_j ^2 \mathcal{F}_j \, \, \forall \, j\neq i.$$
Hence the last equation  of the system is equivalent to
 $$ - \Big( \sum_{k \neq i} \L_k ^2 \mathcal{F}_k\Big) \frac{\xi_i}{\L_i} + \sum_{ j \geq 2; j \neq i}  \L_j ^2 \mathcal{F}_j \frac{\xi_j}{\L_j} = 0  \, \,  \forall \, \, i\geq 2\Longleftrightarrow  \sum_{ j \geq 2}  \L_j ^2 \mathcal{F}_j \frac{\xi_j}{\L_j} = \Big( \sum_{k =1}^m  \L_k ^2 \mathcal{F}_k\Big) \frac{\xi_i}{\L_i} \, \,  \forall \, \, i\geq 2.$$
 Thus we derive that $  \frac{\xi_i}{\L_i} =  \frac{\xi_2}{\L_2} $ for each $i \geq 2$. Putting this information in the last equation implies
 $$ \L_1 ^2 \mathcal{F}_1 \frac{\xi_2}{\L_2} = 0$$
 which implies that $\xi_2=0$ and therefore $H = 0$. This implies that $D \psi_2(Y)$ is invertible and the claim follows from the inverse function theorem.
 \medskip

%%%%%%%%%%%%%%%%%%%%%%%
Next we observe that, using the  coordinates provided by the diffeomorphism $\psi_2(Y)$, the functional  $J_{8m\pi}$ reads as follows
\begin{align*}
 J_{\varrho}(\sum_{i=1}^m\alpha_i\varphi_{a_i, \l_i}\, + \, \ov{w} )   =  & -8 \pi m (1+\ln(m\pi))-8\pi \mathcal{F}_m^K(\ov{a}_1,...,\ov{a}_m) + \sum_{i=1}^m \b_i^2 \\
 & - \sum_{i=2}^m x_i^2 - \Big( \sum_{i=2}^m x_i \Big)^2 +  (4\pi -\s) \frac{ (-\mathcal{L} (Q)) }{ \mathcal{F}^{Q} _1 (q_1)} \frac{\ln {\L}_1}{\L_1 ^2} .
\end{align*}

Finally, it follows from elementary linear algebra, that there exists a change of variables: $ (x_2,\cdots,x_m) \mapsto (y_2, \cdots , y_m)$ such  that $J_{8m\pi}$ reads as
$$  J_{\varrho}(\sum_{i=1}^m\alpha_i\varphi_{a_i, \l_i}\, + \, \ov{w} )   =   -8 \pi m (1+\ln(m\pi))-8\pi \mathcal{F}_m^K(\ov{A}) + \sum_{i=1}^m \b_i^2  - \sum_{i=2}^m y_i^2 +  (4\pi -\s) \frac{ (-\mathcal{L} (Q)) }{ \mathcal{F}^{Q} _1 (q_1)} \frac{\ln {\L}_1}{\L_1 ^2} $$
where $\ov{A}:= (\ov{a}_1,...,\ov{a}_m)$.  This ends the proof.
%%%%%%%%%%%%%%%%%%%%%%%%%
\end{pf}

Next as consequence of Proposition \ref{c:cpatinfinity} of  and the Morse reduction in Lemma \ref{eq:morselem} one derives  the topological contribution of the \emph{critical points at infinity} to the difference of topology between the level sets of the functional $J_\varrho $. Namely we have the following corollary

\begin{cor}\label{l:difftop}
Let $(Q)_{\infty}$ be a critical point at infinity  corresponding to the critical point $Q$ of  $\mathcal{F}^K_m$ at the level  $C_{\infty}(Q)_{\infty}$ with index  $\i_{\infty}(Q)$.  Assume that at the level $C_{\infty}(Q)_{\infty}$ there is no other critical point/critical point at infinity. Then for $\theta$ a small positive number and  a field $\mathbb{F}$, we have that

\begin{align}
 H_l(J_{\varrho}^{C_{\infty}(Q)_{\infty} + \theta},   J_{\varrho}^{C_{\infty}(Q)_{\infty} - \theta }; \mathbb{F} ) =
\begin{cases}
   \mathbb{F} & \mbox{ if  } \quad  l =  \i_{\infty}(Q), \\
  0, & \mbox{otherwise}.
\end{cases}
\end{align}
where $H_l$ denotes the $l-$dimensional homology group with coefficient in the field $\mathbb{F}$.
\end{cor}

%%%%%%%%%%%%%%%%%%%%%%%%%%%%%%%%%%%%%

\section{ Proof of the main results}

%%%%%%%%%%%%%%%%%%%%%%%%%%%%%%%%%%%%%

\begin{pfn}{\bf  of Theorem  \ref{th:crit1}.}
We start by proving the statement $1)$. To this aim   we consider the following sup-approximation of $(MF)$
\begin{equation}\label{eq:TSepsillon}
(MF)_{\mu} \quad  -\D_g u \, = \, 8m\pi (1 + \mu )( \frac{ K e^{u}}{\int_{\Sig} K e^{u}  dV_g} \, - \, 1) \quad \mbox{ in } \Sigma ,
\end{equation}
where $\mu > 0$ is a small positive real number.\\
Regarding Problem $(MF)_{\mu}$ we prove the following : \\
{ \textbf{Claim 1}:  \it For a sequence of $\mu_k \to 0$ , Problem $(MF)_{\mu_k}$ admits a solution $u_{\mu_k}$ whose  generalized  Morse index $Morse(u_{\mu_k})$ is $ 3m$.}

Indeed let $J_{8\pi m(1 +\mu)}$ be the Euler Lagrange functional associated to $(MF)_{\mu}$ then it follows from  \cite{dj, andrea, DM-annals}  that, for large $L$ the sublevel set $ J_{8\pi m(1+{\mu})}^{-L}$ has the same homotopy type as the set of formal barycenters  $ B_m(\Sigma)$ of order $m$ . Note that (see \cite{kk, andrea}), we have
$$ H_{3m -1}( J_{8\pi m(1+\mu)}^{-L}, \mathbb{Z}_2)= H_{3m -1}(B_m(\Sigma), \mathbb{Z}_2)\, \neq \, 0.$$
Let $\varrho_\mu := 8\pi m(1+\mu)$, using the fact that  $J^L_{\varrho_\mu }$ is a contractible set and the exact sequence of the pair $(J^L_{\varrho_\mu}, J^{-L}_{\varrho_\mu} )$ we derive that
$$\begin{CD} 0= H_{3m}(J_{\varrho_\mu }^{L})
\to H_{3m }(J_{\varrho_\mu }^{L}, J_{\varrho_\mu }^{-L}) \to H_{3m-1}(J_{\varrho_\mu }^{-L}) \to H_{3m -1}(J_{\varrho_\mu }^{L})=0.
\end{CD}$$
Hence it follows that
$$ H_{3m}(J^L_{8\pi m(1+\mu)}, J^{-L}_{8\pi m(1+\mu)}) \, \neq 0. $$
Therefore $J_{8\pi m(1+\mu)}$ has a critical point whose Morse index is $ 3m$.

To conclude the proof of the theorem we prove the following claim: \\
 \textbf{Claim 2:} $ \displaystyle{u_{\mu_k}  \to u_{\infty} \mbox{ in } C^{2,\a}(\Sigma),\, }$    where  $u_{\infty} $ is a solution  of Equation $(MF)$.\\
To prove the claim it is enough to rule out the blow up of $u_{\mu_k}$.  Arguing by contradiction it follows from Proposition \ref{p:blowup} that  $ u_{\mu_k} \in V(m,\e)$ for $k$ large . Hence this function has to be written as $$u_{\mu_k} = \sum_{i=1}^m \a_i \varphi_{a_i, \l_i} + w,$$
where the function $w\in E^m_{{A}, \Lambda}$ and satisfies
\begin{equation}\label{****1} \| w \| \leq c \Big(\sum  | \a_k - 1 | + \frac{1}{\l}\Big).\end{equation}
Now, using the fact that $\nabla J_{8\pi m(1+\mu)} (u_{\mu_k})=0$ and Corollary \ref{c:al}, we get
\begin{equation} \label{*1} \sum | \a_i - 1 | \leq c \Big( \frac{1}{\l\ln\l} +  \frac{1}{\ln\l} \sum | \tau_j - \mu + \mu \tau_j | \Big) . \end{equation}
Now, using Proposition \ref{derivelambda} we derive that
\begin{equation} \label{*2}  16 \pi \a_i (\tau_i - \mu + \mu \tau_i ) - 64\pi^2 \sum \frac{\ln\l_j}{\l_j^2} = O\Big( \sum| \a_j - 1 | ^2 \Big) + o\Big(  \frac{\ln\l}{\l^2}  \Big).\end{equation}
Using \eqref{*1} and \eqref{*2}, we get that
\begin{equation}\label{t01} \sum | \tau_i - \mu + \mu \tau_i |  =  O\Big(  \frac{\ln\l}{\l^2}  \Big) . \end{equation}
Hence, using \eqref{*1}, \eqref{t01} and summing \eqref{*2} for $i=1,\cdots,m$, we derive that
$$\sum_{i=1}^m \tau_i - m \mu + \mu  \sum_{i=1}^m \tau_i = 4 \pi m \sum \frac{\ln\l_j}{\l_j^2} + o\Big(  \frac{\ln\l}{\l^2}  \Big).$$
Now using Lemma \ref{sumtau} we obtain
\begin{equation}\label{mu1}\mu = \frac{\pi}{2} \, \mathcal{L}({A }) \frac{\ln \l_1}{\int_\Sigma K e^u } (1+o(1)),\end{equation}
 which implies that $ \mathcal{L}({A }) $ has to be positive.
 Furthermore it follows from Proposition \ref{deriva}, \eqref{*1}, \eqref{t01} and \eqref{****1} that  ${A } := (a_1, \cdots,a_m)$ converges to a critical point of $\mathcal{F}_m^K$.\\
Next expanding the functional $J_{8m\pi(1+\mu_k)}$ in ${V}(m,\e)$ (following the proof of Proposition \ref{devJ}), we obtain
\begin{equation*}
\begin{split}
J_{8m\pi(1+\mu_k)}( u_{\mu_k}) \, & = \, C(m,\mu_k) \,  - \, 8\pi(1 + \mu_k) \mathcal{F}^{K}_{m}(A) \, \, + 16 \pi   \sum_{i=1}^m(\alpha_i-1)^2 \ln \l_i  \\
&  -  \, 4 \pi \,  \sum_{i=1}^m  \wtilde{\tau}_i^2 \, - \, 4\pi \frac{\mathcal{L}(A)}{\mathcal{F}^A_1(a_1)} \frac{\ln \l_1}{\l_1^2 } \, + \, o(\frac{\ln \l}{\l^2}) .
\end{split}
\end{equation*}
Hence  we derive from this expansion, arguing as in Corollary \ref{c:cpatinfinity} that
$$ Morse(u_{\mu_k}) \, = 3 m  - Morse(\mathcal{F}^{K}_m, A). $$
Since $Morse (u_{\e}) = 3m$, we have that  $A$ is a local minimum of  $\mathcal{F}_m^K$  satisfying that $ \mathcal{L}(A) > 0.$ Hence we reach a contradiction to the assumption of Theorem  \ref{th:crit1}.\\
The proof of statement $2)$ follows the same argument as above. The only difference is that we use a sub-approximation:
\begin{equation}\label{eq:TSepsillon-}
(MF)_{\mu} \quad  -\D_g u \, = \, 8m\pi (1 - \mu )\Big( \frac{ K e^{u}}{\int_{\Sig} K e^{u}  dV_g} \, - \, 1\Big) \quad \mbox{ in } \Sigma ,
\end{equation}
where $\mu > 0$ is a small positive real number.
Using the fact that $  H_{3m -4}(B_{m-1}(\Sigma), \mathbb{Z}_2) \neq 0$ (see \cite{kk}, Lemma 8.7) we have that
$$ H_{3m -4}( J_{8\pi m(1-\mu)}^{-L}, \mathbb{Z}_2)= H_{3m -4}(B_{m-1}(\Sigma), \mathbb{Z}_2)\, \neq \, 0.$$
Hence we have that
$$ H_{3m-3}(J^L_{8\pi m(1-\mu)}, J^{-L}_{8\pi m(1-\mu)}) \, \neq 0. $$
Therefore $J_{8\pi m(1-\mu)}$ has a critical point whose Morse index is $ 3m-3$. Next we claim\\
 \textbf{Claim 3:}
$\displaystyle{  u_{\mu_k}  \to u_{\infty} \mbox{ in } C^{2,\a}(\Sigma), \, }$  where $ u_{\infty}$  is a solution  of Equation $(MF)$.\\
By elliptic regularity,  it is enough to rule out the blow up of $u_{\mu_k}$.  Arguing by contradiction we have by Proposition \ref{p:blowup} that for $k$ large
 $ u_{\mu_k} \in V(m,\e)$. It follows then from Proposition \ref{c:cpatinfinity} that $ u_{\mu_k} \in V(m,Q, \e)$, where $Q$ is a critical point of $\mathcal{F}^K_m$ with $\mathcal{L}(Q)< 0.$ Moreover we have that
$$
Morse(u_{\mu_k}) \, = \, 3m - 1- Morse(\mathcal{F}^K_m,Q).
$$
That is $Q$ is a critical point of $\mathcal{F}^K_m$ whose Morse index is 2 and $\mathcal{L}(Q)< 0.$ A contradiction to  the assumption  $2)$ of the theorem. Hence the proof of statement 2)  is complete.
\end{pfn}

%%%%%%%%%%%%%%%%%%%%%%%%%%%%%%%%%%%%%%%%%%%%%%%%%%%%%%%%%%%%%%%

\begin{pfn}{ \bf of Theorem \ref{th:supc3}}
We first observe that since the function $K$ satisfies the non degeneracy condition $(\mathcal{N}_m)$, $J_{\varrho}$  is a Morse function. Moreover the Morse indices of its critical points are uniformly bounded, say by $\bar m$ and it follows from Corollary   \ref{c:cpatinfinity} that the Morse indices of the critical points at Infinity of $J_{8\pi m}$ are bounded by $3m -1$. \\
Without loss of generality, we may assume that $J_{\varrho}$ separates its critical as well as its critical points at Infinity. That is at any critical value there is only one critical point or one critical point at infinity. This can be arranged by perturbing $J_{\varrho}$ slightly in disjoint neighborhoods of its critical points (resp. its critical points at Infinity).
Next we choose $L$ such that all critical values, respectively, critical values at infinity are contained in the open interval $(-L, L)$ and   we order these critical values as
$$ -L < C_1 < \cdots < C_p < L. $$
Now choose regular values $a_0 < \cdots < a_p$ such that
$$
a_0= -L, a_p=L \mbox{ and } \, \, a_{i-1} < C_i < a_i, \, \forall i=1, \cdots,p.
$$
Moreover we  denote by $N_i$ the Morse index of the  critical point $u_i$ such that $J_{\varrho}(u_i) = C_i$, resp. by $N^{\infty}_i$ the  index $\iota_{\infty}$ of the  critical point at infinity $u^{\infty}_i$ such that $C_\infty(u^{\infty}_i) = C_i$ and set $M_i := J_{\varrho}^{a_i}:=\{u: J_{\varrho}(u) < a_i\}$. We notice that it follows from the standard Morse Lemma in the case of usual critical points and from Lemma \ref{eq:morselem} in the case of critical points at infinity,  that $M_i$ is obtained from $M_{i-1}$ by attaching a $N(i)-$cell (resp. $N^{\infty}(i)-$cell).
We set
$$ \theta(i,j):= dim \, H_i(M_j,M_{j-1}) \quad ; \quad \mu(i,j) := dim \,  H_i(M_j,J_{\varrho}^{-L}) $$
and observe that $ \theta(i,j) =  \d_{i,N_j}$ (resp. $ \theta(i,j) =  \d_{i,N^{\infty}_j}$) and $ \mu(i,j) =  0$ if $i > \ov{N}:= \max (\bar m, 3m -1).$
Next, considering the triple $ (M_j, M_{j-1}, J_{\varrho}^{-L})$ we have the exact sequence
$$\begin{CD} 0@>>> H_{*}(M_{j-1}, J_{\varrho}^{-L}) @>>> H_{* }(M_{j}, J_{\varrho}^{-L}),  @>>> H_{*}(M_j,M_{j-1}) \\
 @>>> H_{*-1}(M_{j-1}, J_{\varrho}^{-L})@>>>\cdots   @>>>H_0(M_j,M_{j-1})   \end{CD}$$
We recall that exactness implies the vanishing of corresponding alternating sum of dimensions of the vector spaces in the sequence.

Denoting by $\mathcal{N}_{q,j}$ the kernel of $ H_q(M_{ j-1}, J_{\varrho}^{-L}) \to H_q(M_j,J_{\varrho}^{-L} )$  and $\nu_{q,j} := dim  \mathcal{N}_{q,j}$ and using the exactness of the triple $(M_j; M_{j-1}, J_{\varrho}^{-L})$ we derive that (by grouping 3 terms at a time):
$$ \nu_{q,j} \, =\, \sum_{i=0}^{ q} (-1)^{i + q} \Big( \mu(i,j-1) \, - \mu(i,j)\, + \, \theta(i,j)   \Big ). $$
Summing over $j = 1, \cdots,p$ yields
$$ \sum_{j=1}^{p } \nu_{q,j} \, = \, \sum_{i=0}^{ q} (-1)^{i + q} \Big(\mu(i,0) \, - \mu(i,p)\, + \, \sum_{j=1}^p \theta(i,j)   \Big). $$
Note that $\mu(i,0)=0$ and $\mu(i,p)= \mbox{dim}(H_i(J_\varrho^L,J_\varrho^{-L}))$.
Using the exact sequence of the pair $ (J_{\varrho}^L, J_{\varrho}^{-L}) $ we derive that
\begin{equation} \label{homologie1} \begin{cases} & H_0(J_{\varrho}^{L}, J_{\varrho}^{-L}) \, \simeq \, H_1(J_{\varrho}^{L}, J_{\varrho}^{-L})  \simeq 0  \quad \mbox{ and } \\
& H_i(J_{\varrho}^{L}, J_{\varrho}^{-L}) \, \simeq  \, H_{i-1} (J_{\varrho}^{-L}) \simeq H_{i-1} ( B_{m-1}(\Sigma)) := \beta_{i-1}^{m-1}\quad \forall i \geq 2. \end{cases}\end{equation}
Moreover, denoting by $\nu_i$ the number of critical points of Morse index $i$ resp. by  $\nu^{\infty}_i$ the number of critical points at Infinity of index $i$, it follows that
$$ \sum_{j=1}^{p} \theta(i,j) \, = \nu_i +  \nu^{\infty}_i .$$
Hence we get
\begin{equation} \label{azert}\sum_{j=1}^{p } \nu_{q,j} + \sum_{i=2}^{q} (-1)^{q + i} \beta_{i-1}^{m-1} \, = \, \sum_{( A \in \mathcal{K}^{-}_m; \i_{\infty}(A) \leq q )} (-1)^{\i_{\infty} (A) + q} \, + \, \sum_{i=0}^{q} (-1)^{i+ q} \nu_i. \end{equation}
Now, for $2 \leq k \leq 3m-1$, summing \eqref{azert} for $q=k$ and $q=k-1$, we obtain
$$ \nu_k \, + \, \nu^{\infty}_k  \, - \, \beta^{m-1}_{k-1} \,  =  \sum_{j=1}^{\ov{N}} \nu_{k,j} + \sum_{j=1}^{\ov{N}} \nu_{k-1,j} \geq \, 0, $$
and the statement $(a)$ is proved. \\
The second claim follows by taking $q=\ov{N}$ in \eqref{azert} and using  $\nu_{\ov{N},j}=0$ for each $j$ and
$$ \sum  (-1)^{i}\beta_i^{m-1} = \chi(B_{m-1}(\Sigma)) = 1 - \binom{m - 1 - \chi(\Sigma)}{m-1} .$$
The proof of Theorem \ref{th:supc3} is complete.
\end{pfn}

%%%%%%%%%%%%%%%%%%%%%%%%%

\begin{pfn}{ \bf of Theorem \ref{th:supc4}}
Since there are no critical points at infinity of index $q_0$, it follows from \eqref{azertnew} that:
$$ \nu_{q_0} \, \geq \, \beta^{m-1}_{q_0-1}. $$
Since by assumption $\beta_{q_0-1}^{m-1} \neq 0$, we deduce that Equation $(MF)$ has at least $\beta_{q_0-1}^{m-1}$ solutions. In particular since $\beta^{m-1}_{3m-4} = H_{3m-4}(B_{m-1}(\Sigma),\mathbb{Z}_2) \, \neq 0$, if there are no critical points at infinity of index $3m-3$, we obtain a solution of $(MF)$ whose Morse index is $3m-3$.
\end{pfn}

%%%%%%%%%%%

\begin{pfn}{ \bf of Theorem \ref{th:supc5}}
First observe that if  there is no solution and there is no critical point at infinity of index $q_0$, then there holds:
 \begin{equation}\label{qo} \nu_{q_0,j} =0 \quad \mbox{ and } \quad \nu_{q_0-1,j} =0 \quad \forall \, j.\end{equation}
Indeed, recall that $\nu_{q,j}:= dim\, Ker (i_{q,j})$ where $i_{q,j}$ denotes the map
$$i_{q,j}: \, H_q(M_{j-1},J_\varrho^{-L}) \to H_q(M_{j},J_\varrho^{-L}).$$
The claim is immediate for $q_0$ since there is no critical point or critical point at infinity of index $q_0$ which implies that  $H_{q_0}(M_{j-1},J_\varrho^{-L}) = H_{q_0}(M_{j},J_\varrho^{-L})=0$.
For $q_0-1$, using the exact sequence of the triple $(M_j, M_{j-1}, J_\varrho^{-L})$, we get (since we assumed that there is no critical point/critical point at infinity of index $q_0$)
$$ 0= H_{q_0} (M_j, M_{j-1}) \to H_{q_0-1} (M_{j-1},J_\varrho^{-L}) \to H_{q_0-1}(M_{j},J_\varrho^{-L}), \quad \mbox{ for each }\, j,$$
 which gives the claim for $q_0-1$. \\
Next to prove Theorem \ref{th:supc5}, we  assume that there is no solution and observe that $\beta_i^{m-1} =0$ for each $i > 3m-4$. Furthermore, since there is no critical point at infinity of index $3m$, we get from the above claim \eqref{qo} that $\nu_{3m-1, j} = 0 $ for each $j$.
Hence, summing \eqref{azert} for $q = 3m-1$ and $q=3m-3$, we get
$$0 \leq \sum_{j=1}^{p} \nu_{3m-3,j}  = - (-1)^{3m-3}  \sum_{A\in \mathcal{K}_m^-; \i_\infty(A) = 3m-1\,  et \, \i_\infty(A) = 3m-2}(-1)^{\i_\infty(A)} = \nu _{3m-2}- \nu _{3m-1}.$$

In the same way, summing \eqref{azert} for $q = 3m-1$ and $q=3m-4$, we get
$$\sum_{j=1}^{p} \nu_{3m-4,j} + \beta_{3m-4}^{m-1}   = - (-1)^{3m-4}  \sum_{A\in \mathcal{K}_m^-;  3m-3 \leq  \i_\infty(A)  \leq 3m-1}(-1)^{\i_\infty(A)} = \nu _{3m-3}- \nu _{3m-2}+\nu _{3m-1}.$$
Hence the result follows.
\end{pfn}

\begin{pfn}{ \bf of Theorem \ref{th:new}}
By contradiction, assume that $(MF)$ does not have solutions and observe that, as in \eqref{qo} it follows from assumption  $2)$ that
$$
\forall j, \quad \nu_{q_0 +1,j} \, = \, \nu_{q_0 ,j} \, = \, \nu_{q_0 -1,j} \, = \,  \nu_{q_0 -2,j} = 0.
$$
Hence summing \eqref{azert} for $k=q_0$ and  $k=q_0-1$ we obtain
$ \nu^{\infty}_{q_0} \, = \, \beta^{m-1}_{q_0 -1}.$  A contradiction to  the assumptions $1)$ and $3)$.
\end{pfn}

%%%%%%%%%%%%%%%%%%%%%%%%%%%%%%%%%%%%%
\section{ Appendix }
%%%%%%%%%%%%%%%%%%%%%%%%%%%%%%%%%%%%%

 In this section we collect some technical Lemmas used in this paper.

\begin{lem}\label{AA1}
Let $\varphi_{a,\l}$ and $\d_{a,\l}$ be defined in \eqref{eq:phial} and \eqref{eq:delta}. The following  expansions hold pointwise
\begin{align*}
& \varphi_{a,\l} =  \d_{a,\l} \, + \, \ln\big(\frac{\l^2 }{8}\big) \, + \, 8\pi \, H(a, .) \, + \, 4\pi \frac{\ln\l }{\l^2} \, + \, O\big( \frac{ 1 }{\l^2}\big)\quad \mbox{ in }  \Sigma, \\
& \l \frac{\partial \varphi_{a,\l}}{\partial \l } (x) =  \frac{4}{ 1+ \l^2 \psi_a ^2(x)} - 8\pi \frac{\ln \l }{\l^2} + O\big( \frac{1}{\l^2}\big) \quad \mbox{ for each } x\in \Sigma, \\
&  \frac{1}{\l} \frac{\partial \varphi_{a,\l}}{\partial a } (x) = \frac{1}{\l} \frac{\partial \d_{a,\l}}{\partial a } (x) + 8\pi \frac{1}{\l} \frac{\partial H(a,x)}{\partial a } (x) + O\big( \frac{\ln \l}{\l^2}\big) \quad \mbox{ for each } x\in \Sigma.
\end{align*}
In particular, we have
\begin{align*}
& \varphi_{a,\l} = 8\pi \, G(a, .) \, + \, 4\pi \frac{\ln\l }{\l^2}  \, + \, O\big( \frac{ 1 }{\l^2}\big) \quad  \mbox{ in } \Sigma \setminus B_a(\eta),\\
& \l \frac{\partial \varphi_{a,\l}}{\partial \l } =  - 8\pi \frac{\ln \l }{\l^2} + O\big( \frac{1}{\l^2}\big) \quad \mbox{ in }  \Sigma \setminus B_a(\eta),\\
&  \frac{1}{\l} \frac{\partial \varphi_{a,\l}}{\partial a } =  8\pi \frac{1}{\l} \frac{\partial G(a,.)}{\partial a }  + O\big( \frac{\ln \l}{\l^2}\big) \quad \mbox{ in  } \Sigma  \setminus B_a(\eta).
\end{align*}
\end{lem}
\begin{pf} We remark that the second part of this lemma follows immediately from the first one. Now, we start by proving the first claim.
First, observe that
\begin{equation} \label{edelta}
\int_\Sigma e^{\d_{a,\l} + u_a} dV_g = \int_{B_a(\eta)} e^{\d_{a,\l}} dV_{g_a} + \int_{\Sigma\setminus B_a(\eta)} O\big(\frac{1}{\l^2}\big) = 8\pi + O\big(\frac{1}{\l^2}\big).\end{equation}
Now, let us define the following function $k_{a,\l} := \varphi_{a,\l} - \d_{a,\l} \, - \, \ln\big(\frac{\l^2 }{8}\big) \, - \, 8\pi \, H(a, .)$.
We recall that $\d_{a,\l}$ is a constant function in $ \Sigma \setminus B_a(2\eta)$. Thus, using \eqref{eq:H}, we obtain
\begin{equation*}
 -\D_g k_{a,\l} =  \begin{cases}
& e^{\d_{a,\l} + u_a} -\, \int_{\Sigma}  e^{\d_{a,\l} + u_a} dV_g  +8\pi = O\big(\frac{1}{\l^2}\big) \quad   (\mbox{in } \Sigma \setminus B_a(2\eta)).\\
&  e^{\d_{a,\l} + u_a} -\, \int_{\Sigma}  e^{\d_{a,\l} + u_a} dV_g + e^{u_a} \D_{g_a} \d_{a,\l}+ 8\pi  = O\big(\frac{1}{\l^2}\big) \quad (\mbox{in } B_a(\eta)).\end{cases}
 \end{equation*}
It remains the case of $x\in  B_a(2\eta)\setminus B_a(\eta)$. Using again \eqref{eq:H} we get
\begin{align*}
  -\D_g k_{a,\l} & = O\big(\frac{1}{\l^2}\big) - 8\pi - 2 e^{u_a} \D_{g_a} \ln (1 + \l^2 \psi_a^2) + (8\pi + 2 e^{u_a} \D_{g_a} \ln (\psi_a^2) \\
  & = - 2 e^{u_a} \D_{g_a} \ln \big( \l^2 + \frac{1}{\psi_a^2}\big) + O\big(\frac{1}{\l^2}\big) = O\big(\frac{1}{\l^2}\big) \quad (\mbox{in } B_a(2\eta)\setminus B_a(\eta)).
 \end{align*}
Hence we obtain that $\D_g k_{a,\l} = O(1/\l^2)$ in $\Sigma$ and therefore we get that
$$ k_{a,\l} (x) - \int_{\Sigma} k_{a,\l} (y) dV_g = O\big(\frac{1}{\l^2}\big) \quad \mbox{in }  \Sigma.$$
It remains to estimate the previous integral. Using \eqref{eq:Green}, \eqref{eq:GH} and \eqref{eq:phial}, we get
\begin{align*}
\int_\Sigma k_{a,\l} dV_g & = - \int_\Sigma \d_{a,\l} \, +\, \ln\big(\frac{\l^2 }{8}\big) \, + \, 8\pi \, H(a, .) dV_g \, = \,  2 \int_\Sigma \ln \big(1+  \frac{1}{\l^2 \psi_a^2}\big) dV_g \\
& = 2 \int_{B_a(\eta)} \ln \big( 1+ \frac{ 1 }{\l^2 \psi_a^2}\big) dV_{g_a}  + O\big(  \int_{B_a(\eta)} \ln \big( 1+ \frac{ 1 }{\l^2 \psi_a^2}\big) |e^{-u_a} - 1 | dV_{g_a} + \frac{1}{\l^2}\big)\\
& = \frac{ 4 \pi }{ \l^2} \int_0^{\l \eta}\ln \big( 1 + \frac{1}{r^2}\big) r dr + O\big(\frac{1}{\l^2}\big) =   4 \pi\frac{ \ln \l}{ \l^2}  + O\big(\frac{1}{\l^2}\big) .
\end{align*}
 Hence the proof of the first claim follows.\\
 The other claims can be proved in the same way.
 \end{pf}

\begin{lem}\label{AA2}
Let $\varphi_{a,\l}$  be defined in \eqref{eq:phial}.  The following  expansions hold:
$$ \| \varphi_{a,\l} \| ^2 = 32 \pi \ln \l + 64 \pi^2 H(a,a) -16 \pi + 64 \pi^2 \frac{\ln \l}{\l^2} + O\big( \frac{1}{\l^2}\big) ,$$
$$ \langle \varphi_{a,\l},  \l \frac{\partial \varphi_{a,\l}}{\partial \l } \rangle_{g} = 16 \pi  - 64 \pi^2 \frac{\ln \l}{\l^2} + O\big( \frac{1}{\l^2}\big) ,$$
$$ \langle \varphi_{a_j,\l_j},   \varphi_{a_i,\l_i} \rangle_{g} =  64\pi^2 \, G(a_j, a_i) \, + \, 32\pi^2 \frac{\ln\l_j}{\l_j^2} \, + \, 32\pi^2 \frac{\ln\l_i}{\l_i^2}  \, + \, O\big( \frac{ 1 }{\l_j^2} +  \frac{1}{\l_i^2} \big),$$
$$ \langle \varphi_{a_j,\l_j},  \l_i \frac{\partial \varphi_{a_i,\l_i}}{\partial \l_i } \rangle_{g} = - 64 \pi^2 \frac{\ln \l_i}{\l_i^2} + O\big( \frac{1}{\l_i^2} +  \frac{1}{\l_j^2} \big) .$$
\end{lem}
\begin{pf} Since $ \int_\Sigma  \varphi_{a,\l} dV_g = 0$, using Lemma \ref{AA1} and \eqref{edelta}, we get
\begin{align*}
 \| \varphi_{a,\l} \| ^2 &  =\int_\Sigma e^{u_a} e^{\d_{a,\l}} \big( \d_{a,\l} \, + \, \ln\big(\frac{\l^2 }{8}\big) \, + \, 8\pi \, H(a, .) \, + \, 4\pi \frac{\ln\l }{\l^2} \, + \, O\big( \frac{ 1 }{\l^2}\big)\big)dV_g  \\
& = 32\pi^2 \frac{\ln\l }{\l^2} \, + \, O\big( \frac{ 1 }{\l^2}\big)+ 2 \int_\Sigma e^{u_a} e^{\d_{a,\l}}  \ln \big( \frac{\l^2\psi_a^2}{1+ \l^2\psi_a^2}\big)dV_g + 8\pi \int_\Sigma e^{u_a} e^{\d_{a,\l}} G(a,.)dV_g
\end{align*}
Note that (using the fact that $\int_\Sigma G(a,.) dV_g = 0$)
$$ 8\pi \int_\Sigma e^{u_a} e^{\d_{a,\l}} G(a,.)dV_g = 8\pi \varphi_{a,\l}(a) = 32 \pi \ln \l + 64 \pi^2 H(a,a)  + 32 \pi^2 \frac{\ln \l}{\l^2} + O\big( \frac{1}{\l^2}\big).$$
Furthermore, we have
\begin{align*}
\int_\Sigma e^{u_a} e^{\d_{a,\l}}  \ln \big( \frac{\l^2\psi_a^2}{1+ \l^2\psi_a^2}\big)dV_g & = \int_{B_a(\eta)}  e^{\d_{a,\l}} \ln \big( \frac{\l^2 | x-a|^2}{1+ \l^2 | x-a |^2}\big)dV_{g_a}  + O\big( \frac{1}{\l^4}\big)   \\
& = \int_{B(0,\eta)} \frac{8 \l^2}{ (1+ \l^2 | x | ^2)^2}\ln \big( \frac{\l^2 | x |^2}{1+ \l^2 | x |^2}\big) dx + O\big( \frac{1}{\l^4}\big)   \\
& = - 8\pi + O\big( \frac{1}{\l^2}\big)  .
\end{align*}
Hence the proof of the first claim follows.\\
Concerning the second one, using the fact that  $ \int_\Sigma  \l \frac{\partial \varphi_{a,\l}}{\partial \l } dV_g = 0$, we have
\begin{align*}
 \langle \varphi_{a,\l},  \l \frac{\partial \varphi_{a,\l}}{\partial \l } \rangle_{g} &  = \int_\Sigma e^{u_a} e^{\d_{a,\l}} \bigg( \frac{4}{ 1+ \l^2 \psi_a ^2(x)} - 8\pi \frac{\ln \l }{\l^2} + O\big( \frac{1}{\l^2}\big) \bigg)dV_g \\
& = - 64 \pi^2 \frac{\ln \l }{\l^2} + O\big( \frac{1}{\l^2}\big)  +  \int_\Sigma e^{u_a} e^{\d_{a,\l}}\frac{4}{ 1+ \l^2 \psi_a ^2(x)} dV_g .
\end{align*}
Observe that, we have
$$ \int_\Sigma e^{u_a} e^{\d_{a,\l}}\frac{4}{ 1+ \l^2 \psi_a ^2(x)} dV_g  = \int_{B(0,\eta)}  \frac{ 32 \l^2}{ ( 1+ \l^2 | x |^2 )^3} dx + O\big( \frac{1}{\l^4}\big)  = 16\pi + O\big( \frac{1}{\l^4}\big) .$$
Thus the proof of the second claim follows.
Now, we will focus on the third claim. Using  $\int_\Sigma |  \varphi_{a,\l} | =O(1)$ and Lemma \ref{AA1}, it holds
\begin{align*}
\langle \varphi_{a_j,\l_j},   \varphi_{a_i,\l_i} \rangle_{g} & = \int_{B_{a_j}(\eta)} e^{u_{a_j}} e^{\d_{a_j,\l_j}}  \bigg( 8\pi \, G(a_i, .) \, + \, 4\pi \frac{\ln\l_i }{\l_i^2}  \, + \, O\big( \frac{ 1 }{\l_i^2}\big)\bigg) dV_g+ O\big( \frac{1}{\l_j^2}\big) \\
& = 32\pi^2 \frac{\ln\l_i }{\l_i^2} + O\big( \frac{1}{\l_j^2} + \frac{1}{\l_i^2}\big) + 8\pi  \int_{B_{a_j}(\eta)} e^{u_{a_j}} e^{\d_{a_j,\l_j}}    G(a_i, .) dV_g .
\end{align*}
Observe that (using $\int_\Sigma G(a_i,.) dV_g =0$)
\begin{align*}
\int_{B_{a_j}(\eta)} e^{u_{a_j}} e^{\d_{a_j,\l_j}}    G(a_i, .) dV_g & = \int_{\Sigma} e^{u_{a_j}}  e^{\d_{a_j,\l_j}}    G(a_i, .)  dV_{g} + O\big( \frac{1}{\l_j^2} \big)\\
&  =  \varphi_{a_j,\l_j}(a_i)+ O\big( \frac{1}{\l_j^2} \big) = 8\pi \, G(a_j, a_i) \, + \, 4\pi \frac{\ln\l_j}{\l_j^2}  \, + \, O\big( \frac{ 1 }{\l_j^2} \big).
\end{align*}
Hence the proof of this claim follows.
Concerning the last one, it holds
\begin{align*}
\langle \varphi_{a_j,\l_j},  \l_i \frac{\partial \varphi_{a_i,\l_i}}{\partial \l_i } \rangle_{g} & =  \int_{B_{a_j}(\eta)} e^{u_{a_j}} e^{\d_{a_j,\l_j}}  \bigg( -8\pi \frac{\ln\l_i }{\l_i^2} + O\big( \frac{ 1 }{\l_i^2} \big)\bigg)dV_g + O\big( \frac{ 1 }{\l_j^2} \big)\\
& = -64 \pi^2 \frac{\ln\l_i }{\l_i^2} + O\big( \frac{ 1 }{\l_j^2} + \frac{ 1 }{\l_i^2} \big).
\end{align*}
Thereby the proof of this lemma follows.
\end{pf}

From Lemma \ref{AA1}, it is easy to get the following expansion
\begin{lem}\label{eu} Let ${u}:= \sum_{j=1}^m \a_j \varphi_{a_j, \l_j}$.  In $B_{a_i}(\eta)$, it holds
\begin{align*}
 K e^{u} & =  \frac{ \lambda_i^{4\alpha_i}\mathcal{F}_i^{A}g_i^{A}}{(1+ \lambda_i^2| y_{a_i}(.) | ^2)^{2\alpha_i}}  \Big(1+ 4\pi \sum_{j=1}^m  \frac{\ln\l_j}{\l_j^2} + O\big(\frac{1}{ \lambda^2} + \sum | \a_j -1 | \frac{\ln\l_j}{\l_j^2} \big) \Big)\\
&  =\frac{ \lambda_i^{4\alpha_i}  \mathcal{F}_i^{A}g_i^{A}}{(1+ \lambda_i^2 | y_{a_i}(.) |^2)^{2\alpha_i}}  \Big(1+ 4\pi \sum_{j=1}^m  \frac{\ln\l_j}{\l_j^2} \Big) + O \bigg(\Big( \frac{1}{ \lambda^2} + \sum | \a_j -1 | \frac{\ln\l_j}{\l_j^2}\Big) Ke^{u}\bigg), \end{align*}
where the functions $\mathcal{F}^{A}_i$ (with ${A} = (a_1, \cdots, a_m)$) and $g_i^{A}$  are defined as follows
\begin{equation}  \label{eq:fkmi2}
\left\{ \begin{array}{cc}
  \mathcal{F}^{A}_i(x)  :=  K(x) \exp \Big( 8 \pi H(a_i,x) +8 \pi \sum_{j \neq i} G(a_j,x)\Big), \vspace{1mm} \\
 g_i^{A}(x) := \exp\Big(8\pi(\alpha_i-1) H(a_i,x)+8\pi  \sum_{j \ne i}(\alpha_j-1)G(a_j,x)\Big).  \end{array} \right.
\end{equation}
\end{lem}

\begin{lem} \label{intkeu}
Let $u:=\sum_{i=1}^m \alpha_i \varphi_i  \in V(m,\varepsilon)$. Then,
\begin{eqnarray}
\int_{\Sigma} Ke^u dV_g= \pi \sum_{i=1}^m\frac{  \lambda_i^{4\alpha_i-2} }{2\alpha_i-1} (\mathcal{F}_i^{A}g_i ^{A})(a_i)  +O \Big(1+ \sum_{k=1}^m  \lambda_k^{4\alpha_k-2}  \Big(|\alpha_k-1|^2 + \frac{\ln  \lambda_k}{ \lambda_k^{2}}\Big)\Big).  \label{eq:Raffaella}
\end{eqnarray}

If  $ \sum_{i=1}^m |\alpha_i-1|\ln \lambda_i = o_{\varepsilon} (1)$,  the expansion~\eqref{eq:Raffaella} can be improved as follows
\begin{eqnarray}
\int_{\Sigma} Ke^u dV_g &=& \pi \sum_{i=1}^m  \frac{ \lambda_i^{4\alpha_i-2}}{2\alpha_i-1} ( \mathcal{F}_i^{A} g_i^{A}) (a_i) +  \frac{\pi}{2}\sum_{i=1}^m \Big( \Delta \mathcal{F}_i^{A} (a_i) - 2 K_g(a_i) \mathcal{F}_i^{A} (a_i) \Big) \ln \lambda_i   \nonumber  \\
  && + 4\pi^2 \Big( \sum_{j=1}^m\frac{\ln\l_j}{\l_j^2} \Big) \sum_{i=1}^m   \lambda_i^{4\alpha_i-2} \mathcal{F}_i^{A}(a_i) +  O \Big( 1+ \sum_{k=1}^m  |\alpha_k-1|\ln^2 \lambda_k \Big). \label{eq:Raffaella2}
\end{eqnarray}
\end{lem}

\begin{pf}
Let $B_i:= B_{a_i}(\eta)$. Firstly, note that Lemma~\ref{AA1} shows that $e^u$ is bounded on $\Sigma \setminus   (\cup B_i)$. Hence the integral in this set is bounded. Now, using Lemma \ref{eu}, we have
\begin{align*}
&  \int_{B_i}  Ke^{u} dV_g \\
&  =\int_{B_i} \frac{ \lambda_i^{4\alpha_i}  \mathcal{F}^{A}_i g_i^{A} }{(1+ \lambda_i^2 | y_{a_i} |^2)^{2\alpha_i}} \Big(1+ 4\pi \sum_{j=1}^m  \frac{\ln\l_j}{\l_j^2} \Big) dV_g + O\bigg(\big(\frac{1}{ \lambda^2}   + \sum | \a_j -1 | \frac{\ln\l_j}{\l_j^2} \big) \int_{B_i} Ke^{u}dV_g\bigg)\\
 & = \Big(1+ 4\pi \sum_{j=1}^m  \frac{\ln\l_j}{\l_j^2} \Big) \int_{B_i} \frac{ \lambda_i^{4\alpha_i}  \mathcal{F}^{A}_i g_i^{A}e^{-u_{a_i}}}{(1+ \lambda_i^2 | y_{a_i} (x)|^2)^{2\alpha_i}} dV_{g_{a_i}} + O\bigg(\big(\frac{1}{ \lambda^2}   + \sum | \a_j -1 | \frac{\ln\l_j}{\l_j^2} \big) \int_{\Sigma} Ke^{u}dV_g\bigg).
 \end{align*}
 Finally, we have
\begin{align}\int_{B_i} \frac{ \lambda_i^{4\alpha_i}\mathcal{F}^{A}_i g_i^{A}e^{-u_{a_i}} }{(1+ \lambda_i^2 | y_{a_i} (x)|^2)^{2\alpha_i}} dV_{g_{a_i}}
& =  \lambda_i^{4\alpha_i-2} (\mathcal{F}_i^{A} g_i^{A}e^{-u_{a_i}}) (a_i)\int_{B_i} \frac{ \lambda_i^{2}}{(1+ \lambda_i^2 | y_{a_i} (x)|^2)^{2\alpha_i}} dV_{g_{a_i}} \label{111}\\
  & +  \lambda_i^{4\alpha_i-2} O \left(\int_{B_i} \frac{ \lambda_i^{2}| y_{a_i} (x)|^2}{(1+ \lambda_i^2 | y_{a_i} (x)|^2)^{2\alpha_i}} dV_{g_{a_i}} \right).\nonumber\end{align}
 Note that ${u_{a_i}}(a_i)=0$. Furthermore, easy computations imply
\begin{equation}\label{112} \int_{B(0,\eta)} \frac{ \lambda_i^{2}}{(1+ \lambda_i^2|x|^2)^{2\alpha_i}} dx = \frac{\pi}{2\alpha_i-1}+ O \Big(\frac{1}{ \lambda_i^{4\alpha_i-2}} \Big). \end{equation}
To estimate the other integral, we introduce the following function
 \begin{equation}  \label{e:xi}    \xi(x)= \frac{1}{(1+ \lambda^2|x|^2)^{2(\alpha-1)}}. \end{equation}
If $(\alpha-1)\ln  \lambda$ is small, one obtains that $\xi=1+o(1)$ uniformly on $B(0, \eta)$. In general, we have
\begin{equation} \label{e:nxi} |\xi(x)-1|=|\xi(x)-\xi(0)| \, = \, \Big|  \int_0^1   \frac{4(1-\alpha) \lambda^2t|x|^2}{(1+ \lambda^2t^2|x|^2)^{2\alpha-1}}dt \Big|   \, \leq \,  c|\alpha-1| \sqrt{ \lambda |x|}.\end{equation}
Now, using \eqref{e:nxi}, we obtain
\begin{align*}
\int_{B(0,\eta)} \frac{ \lambda_i^2|x|^2}{(1+ \lambda_i^2|x|^2)^{2\alpha_i}} dx & = \int_{B(0,\eta)} \frac{ \lambda_i^2|x|^2}{(1+ \lambda_i^2|x|^2)^{2}} dx + \int_{B(0,\eta)} \frac{ \lambda_i^2|x|^2}{(1+ \lambda_i^2|x|^2)^{2}}(\xi_i -1) dx \\
 & =O \Big(\frac{\ln  \lambda_i}{ \lambda_i^2} \Big) +     O \Big(\frac{|\alpha_i -1|}{ \lambda_i^{3/2}} \Big).
\end{align*}
Hence, the proof of \eqref{eq:Raffaella} follows.
Now, we will focus on \eqref{eq:Raffaella2}. In this case, we have $| \a_i - 1 | \ln \l_i$ is small for each $i$ and we need to improve the estimate of \eqref{111}. Using \eqref{112}, we obtain
\begin{align*}\int_{B_i} & \frac{ \lambda_i^{4\alpha_i}\mathcal{F}^{A}_i g_i^{A}e^{-u_{a_i}} }{(1+ \lambda_i^2 | y_{a_i} (x)|^2)^{2\alpha_i}}dV_{g_{a_i}} = \frac{\pi}{2\alpha_i-1} \lambda_i^{4\alpha_i-2} (\mathcal{F}_i^{A} g_i^{A}) (a_i)  + O(1)\\
& +  \frac{1}{4} \D_{g_{a_i}} (\mathcal{F}_i^{A} g_i^{A}e^{-u_{a_i}}) (a_i) \int_{B(0,\eta)}  \frac{\l_i^{4\a_i} | x |^2 }{ (1+ \l_i^2 | x | ^2 )^{2\a_i}} dx + O\Big( \int_{B(0,\eta)}  \frac{\l_i^{4\a_i} | x |^3 }{ (1+ \l_i^2 | x | ^2 )^{2\a_i}} dx \Big).
\end{align*}
Observe that
\begin{align*}
& \int_{B(0,\eta)}  \frac{\l_i^{4\a_i} | x |^3 }{ (1+ \l_i^2 | x | ^2 )^{2\a_i}} dx = \l_i^{4\a_i -5} 2\pi \int_0 ^{\l_i \eta}. \frac{r^4}{(1+r^2)^{2\a_i}} dr = O(1),\\
& \int_{B(0,\eta)}  \frac{\l_i^{4\a_i} | x |^2 }{ (1+ \l_i^2 | x | ^2 )^{2\a_i}} dx = 2 \pi \ln\l_i + O( 1+| \a_i -1 | \ln^2\l_i ) \quad (\mbox{by using }  | \a_i - 1 | \ln\l_i \mbox{ is small}).
\end{align*}
To conclude, we need the following information. We know that the function $\wtilde{u}_a(x):=u_a (y_a(x))$ (for $x\in B(a,\eta) \subset \R^2$) satisfies $$ -\D \wtilde{u}_a = -2 K_g(y_a^{-1}(.)) e^{-\wtilde{u}_a} \quad \mbox{ in } B(a,\eta) $$ and therefore we derive that
$$ \D_{g_{a_i}} (\mathcal{F}_i^{A} g_i^{A}e^{-u_{a_i}}) (a_i) = \D_{g_{a_i}} \mathcal{F}_i^{A}  (a_i) -  2 \mathcal{F}_i^{A}(a_i) K_g(a_i) + O(\sum | \a_j -1 | ),$$
by using the fact that
$$g_i^{A} (a_i)= 1 + O\Big( \sum | \a_j - 1 | \Big), \quad | \nabla g_i^{A} (a_i)| =  O\Big( \sum | \a_j - 1 | \Big) , \quad | \D g_i^{A} (a_i)|=  O\Big( \sum | \a_j - 1 | \Big).$$
Summing the above estimates, the result follows.
\end{pf}

\begin{lem}\label{l:Ual}
Let  $A$ and $\L$ satisfying the  balancing condition \eqref{eq:balanc-cond} then  the approximate solution $U_{A,\L}$ defined in \eqref{eq:Ual} satisfies the following equation
$$
-\D_g U_{A,\L} \, + \, 8 \pi m \, = \, 8 \pi m \frac{K e^{U_{A,\L}}}{\int_{\Sig} K e^{U_{A,\L}}} \, + \, f_{A,\L},
\quad \mbox{  where } $$
$$
f_{A.\L} \rightharpoonup 0 \mbox{ weakly in }  H^1(\Sig ) \mbox{ as } \l_i \to \infty \, \forall \, i=1, \cdots,m.
$$
\end{lem}

\begin{pf}
  We first notice that according to Equation \eqref{eq:phial} satisfied by $\varphi_{a,\l}$, using \eqref{edelta}, the function $U_{A,\l}$ satisfies the equation
  \begin{equation}\label{eq:Ual1}
  -\D_g U_{A,\L} \, + \, 8 \pi m \, + \, O( \sum_{i=1}^{m} \frac{1}{\l_i^2}) = \, \sum_{i=1}^{m} e^{\d_{a_i,\l_i }+ u_{a_i}}.
  \end{equation}
Moreover it follows from Lemmas \ref{AA1} and \ref{intkeu} that in $\Sig \setminus \bigcup_{i=1}^m B_{\eta}(a_i)$ we have that
$$
8 \pi m \frac{K e^{U_{A,\L}}}{\int_{\Sig} K e^{U_{A,\L}}} \, = \, O(\frac{\ln \l}{\l^2}) .
$$
Furthermore it follows from  Lemmas \ref{eu} and \ref{intkeu} that in $B_{\eta}(a_i)$ there hold
\begin{equation*}
 K e^{U_{a,\l}} \, = \, \frac{1}{8} \l_i^2 \mathcal{F}_i^A(a_i) e^{\d_{a_i,\l_i}} ( 1 + O(\frac{\ln \l}{\l^2})) \quad \mbox{ and } \quad
 \int_{\Sig} K e^{U_{A,\l}} \, = \, \pi \sum_{i=1}^{m} \l_i^2 \mathcal{F}_i^A(a_i)( 1 + O(\frac{\ln \l}{\l^2})). \end{equation*}
Therefore in $B_{\eta}(a_i)$  we have that
$$
8 \pi m \frac{K e^{U_{A,\L}}}{\int_{\Sig} K e^{U_{A,\L}}}\, = \, \frac{m\l_i^2 \mathcal{F}^A_i(a_i) e^{\d_{a_i,\l_i}}}{\sum_{j=1}^{m} \l_j^2 \mathcal{F}^A_j(a_j)} ( 1 + O(\frac{\ln \l}{\l^2})) .
$$
Hence we derive that
$$ \sum_{j=1}^{m} e^{\d_{a_j,\l_j }+ u_{a_j} } \, - \, 8 \pi m \frac{Ke^{U_{A,\L}}}{\int_{\Sig} K e^{U_{A,\L}}} \, =  e^{\d_{a_i,\l_i}}( \tilde{\tau}_i + O( |x-a_i|^2\frac{\ln \l}{\l^2})) \quad \mbox{ in } B_\eta(a_i), $$
where ${\tau}'_i$ is defined in Proposition \ref{devJ} by taking $\a_k = 1$ for each $k$.\\
We notice that  it follows from the balancing condition \eqref{eq:balanc-cond} that  $\forall i=1, \cdots,m$ we have that $|{\tau}'_i| = o_{\l}( 1 )$.
Summarizing we have that $U_{A,\L}$ satisfies the equation
\begin{align*}
& -\D_g U_{A,\L} \, + \, 8 \pi m \, = \, 8 \pi m \frac{K e^{U_{A,\L}}}{\int_{\Sig} K e^{U_{A,\L}}} \, + \, f_{A,\L},
\quad \mbox{ where }\\
& f_{A,\L} \, = \,
\begin{cases}
 O({\ln \l}/{\l^2})  \mbox{ in  } \Sig \setminus \bigcup_{i=1}^m B_{\eta}(a_i) \\
   e^{\d_{a_i,\l_i}}( \tilde{\tau}_i + O({ |x-a_i|^2 + \ln \l_i}/{\l_i^2})) \mbox{ in } B_{\eta}(a_i).
\end{cases}
\end{align*}
Hence the Lemma follows.
\end{pf}

Lastly for the sake of completeness  we provide the following characterization of  approximate blowing up solutions of Equation $(MF)$. Namely we prove:
%%%%%%%%%%%%%%%%%%%%%%%%%%%%%%%%%%%%%%%%%%%%%%%%%%%%%%%
\begin{pro}\label{p:blowup}
Let $(\Sig,g)$ be a closed surface of unit volume and  $u_k  $ be  a blowing up solution  of
$$
-\D_g u_k \, = \, \varrho_k( \frac{K e^{u_k}}{\int_{\Sig} K e^{u_k} } \, - 1) \quad \mbox{ in } \Sig ; \qquad \int_{\Sig} u_k dV_g = 0,
$$
where $\varrho_k \to 8 \pi m, \, m \in \N$. Then for $\e$ small and $k$ large we have that $u_k \in V(m,\e)$.
\end{pro}
%%%%%%%%%%%%%%%%%%%%%%%%%%%%%%%%%%%%%%%%%%%%
\begin{pf}
We set $\tilde{u}_k:= u_k  \, - \,  \ln(\int_{\Sig} K e^{u_k})$ and observe that $ \tilde{u}_k$ satisfies the equation

$$
- \D_g \tilde{u}_k \, = \, \varrho_k( Ke^{\tilde{u}_k} -1) \, \mbox{ in } \Sig.
$$
Now it follows from the refined blow up analysis performed in \cite{LS,cl1,cl2} that $\tilde{u}_k$ blows up at $m$-points $a_1, \cdots,a_m \in \Sig$  with comparable blow up rates $\l_i:= e^{u_k(a_i)}$ such that $d_g(a_i,a_j) \geq 2 \eta$ for some $\eta > 0$. Moreover we have that
$$
e^{\tilde{u}_k} \, = \, O( \sum_{i=1}^{m} \frac{1}{\l_i^2}) \, \mbox{ in } \Sig \setminus \bigcup_{i=1}^m B_{\eta}(a_i)
$$
and inside the balls  $B_i:= B_{\eta}(a_i)$ there holds (See Theorem 1.4 in \cite{cl1})
$$
\tilde{u}_k \, + \, \ln (\varrho_k K(a_i)) \, = \, {\d_{a_i,\l_i}} \, + O(d_g(a_i,x)) .
$$
Next we set
$$
\tilde{v}_k := \tilde{u}_k \, - \, \sum_{i=1}^{m} \varphi_{a_i,\l_i}
$$
and we notice that $\tilde{v}_k$ satisfies the equation
\begin{equation}\label{eq:vktilde}
  -\D_g \tilde{v}_k  \, + \, \varrho_k - 8 m \pi \,    = \varrho_k K e^{\tilde{u}_k} - \sum_{i=1}^{m} e^{\d_{a_i,\l_i }+ u_{a_i}} + O(\frac{1}{\l^2})  := f_k.
\end{equation}
Next we claim that: \\
\textbf{Claim}: For $1 < p < 2$ there exists a constant $C > 0$ such that
$$
\|f_k\|_{L^p} \, \leq \, C \sum_{i=1}^{m} \frac{1}{\l_i^{(2/p)- 1}}.
$$
Indeed we have that
$$
\int_{\Sig} |f_k|^p dV_g \, = \, \sum_{i=1}^{m} \int_{B_i} |f_k|^p dV_g  \, + \, O( \sum_{i=1}^{m} \frac{1}{\l_i^{2p}}).
$$
Furthermore  in $B_i$ we have, in  geodesic local coordinates around $a_i$, that
\begin{align}\label{eq:fk}
  f_k & =   e^{\tilde{u}_k + \ln (\varrho_k K)} \,  - \sum_{i=1}^{m} e^{\d_{a_i,\l_i}+u_{a_i}} + O(\frac{1}{\l^2}) \nonumber  \\
   & =  e^{ \d_{a_i,\l_i} + \ln (\varrho_k K) - \ln (\varrho_k K(a_i) + O(|x-a_i|))} - e^{\d_{a_i,\l_i}+u_{a_i}} + O(\frac{1}{\l^2}) \nonumber  \\
   & = O(|x-a_i| e^{\d_{a_i,\l_i}}  + \frac{1}{\l^2}) .
\end{align}
Hence  for $ 1 < p < 2$  and a constant $C > 0$ there holds
$$
\int_{\Sig} |f_k|^p dV_g \, \leq \, C  \, \sum_{i=1}^{m} \frac{1}{\l_i^{2-p}}.
$$
Next setting $v_k:= \tilde{v}_k \, + \, \ln (\int_{\Sigma } K e^{u_k}) \, = \, u_k - \sum_{i=1}^{m} \varphi_{a_i,\l_i}$ we have that $v_k$ satisfies the equation

$$
- \D_g \, v_k+ \, \varrho_k  - 8 \pi m \, = \, f_k \, \mbox{ in } \Sig;  \, \int_{\Sig} v_k dV_g \, = \, 0 ,
$$
where $f_k \in L^p(\Sig) $ for $p \in (1,2)$. Hence it follows from the Calderon-Zygmund  a priori estimate that
$$
\|v_k\|_{W^{2,p}} \, \leq \, C \, \|f_k\|_{L^p} \, \leq \, C \, \sum_{i=1}^{m} \frac{1}{\l_i^{(2/p) -1}}.
$$
Therefore it follows from the Sobolev embedding  $W^{2,p}(\Sig ) \hookrightarrow W^{1,2}(\Sig)$ that
$$
\|v_k \|_{H^1} \, \leq \,  \sum_{i=1}^{m} \frac{1}{\l_i^{(2/p) -1}}.
$$
Hence choosing $p= \frac{4}{3}$ we obtain that
$$
\|u_k \, - \, \sum_{i=1}^{m} \varphi_{a_i,\l_i} \| \, \leq   \, C \, \sum_{i=1}^{m} \frac{1}{\sqrt{\l_i}}.
$$
Therefore the proposition is fully proven.
\end{pf}

%%%%%%%%%%%%%%%%%%%%%%%%%%%%%%%%%%%%%%%%%%%%%%%
%%%%%%%%%%%%%%%%%%%%%%%%%%%%%%%%%%%%%%%%%%%%%%%%

\bigskip
\bigskip

$$
\begin{minipage}[l]{7.5cm}
    {\small Mohameden Ahmedou}\\
    {\footnotesize Mathematisches Institut} \\
     {\footnotesize der Justus-Liebig-Universit\"at Giessen}\\
    {\footnotesize Arndtsrasse 2, D-35392 Giessen}\\
    {\footnotesize Germany}\\
    {\tt\footnotesize Mohameden.Ahmedou@math.uni-giessen.de  }
    {}
\end{minipage}
\quad
\begin{minipage}[l]{7.5cm}
    {\small Mohamed Ben Ayed}\\
    {\footnotesize Universit\'e de Sfax,  Facult\'e des Sciences de Sfax}\\
    {\footnotesize D\'epartement de Math\'ematiques}\\
    {\footnotesize Route de Soukra, Sfax, Tunisia}\\
    {\tt\footnotesize  Mohamed.Benayed@fss.rnu.tn}
    {}
\end{minipage}
$$

\end{document}